\documentclass[%
  onecolumn
   , colorlinks 
]{mpi2015-cscpreprint}

\usepackage[american]{babel}

\usepackage{graphicx}

\usepackage{amssymb}
\usepackage{amsthm}

\newcommand{\ode}{\texttt{ODE}}
\newcommand{\pde}{\texttt{PDE}}
\newcommand{\pinn}{\texttt{PINN}}
\newcommand{\fdm}{\texttt{FDM}}
\newcommand{\fem}{\texttt{FEM}}
\newcommand{\dnn}{\texttt{DNN}}

\newcommand{\sindy}{\texttt{SINDy}}
\newcommand{\propod}{\texttt{POD}}
\newcommand{\nn}{\texttt{NN}}
\newcommand{\deim}{\texttt{DEIM}}
\newcommand{\eim}{\texttt{EIM}}
\newcommand{\qdeim}{\texttt{Q-DEIM}}
\newcommand{\maxm}{\texttt{max}}
\newcommand{\minm}{\texttt{min}}
\newcommand{\SVD}{\texttt{SVD}}

\newcommand{\gspinn}{\texttt{GS-PINN}}

\newcommand{\gnsindy}{\texttt{GN-SINDy}}
\newcommand{\deepymod}{\texttt{DeePyMoD}}
\newcommand{\lasso}{\texttt{LASSO}}
\newcommand{\stridge}{\texttt{STRidge}}
\newcommand{\ols}{\texttt{OLS}}

\usepackage[vlined,longend,linesnumbered,ruled]{algorithm2e}

\usepackage[american]{babel}

\usepackage{graphicx}

\usepackage{amssymb}
\usepackage{amsthm}
\usepackage{amsmath}

\usepackage{caption}
\usepackage{subcaption}

\usepackage{graphics} 
\usepackage{epsfig} 
\usepackage{tikz}
\usepackage{xcolor}
\usepackage[latin1]{inputenc}
\usepackage[T1]{fontenc}

\usepackage{cleveref}

\usepackage{multirow}

\usepackage{comment}

\usepackage{longtable,array,amsmath,booktabs,tabularx,ragged2e}

\usetikzlibrary{arrows,automata, quotes, calc, trees, positioning, chains, shapes.geometric, decorations.pathreplacing, decorations.pathmorphing, shapes, matrix, shapes.symbols}

\usepackage{circuitikz}

\usetikzlibrary{arrows}
\usepackage[noend]{algpseudocode}

\usepackage{lscape}
\usepackage[colorinlistoftodos,prependcaption,textsize=footnotesize]{todonotes}

\usepackage{hyperref}
\usepackage{enumitem}

\usepackage[software, hardware]{mymacros}

\begin{document}
  

\title{GN-SINDy: Greedy Sampling Neural Network in Sparse Identification of Nonlinear Partial Differential Equations}

\author[$\ast$]{Ali Forootani}
\affil[$\ast$]{Max Planck Institute for Dynamics of Complex Technical Systems, 39106 Magdeburg, Germany.\authorcr \email{forootani@mpi-magdeburg.mpg.de}, \orcid{0000-0001-7612-4016}}

\author[$\dagger$]{Harshit Kapadia}
\affil[$\dagger$]{ \email{kapadia@mpi-magdeburg.mpg.de}, \orcid{0000-0003-3214-0713}}

\author[$\ddagger$]{Sridhar Chellappa}
\affil[$\ddagger$]{
	\email{chellappa@mpi-magdeburg.mpg.de}, \orcid{0000-0002-7288-3880}}

\author[$\ddagger\ddagger$]{\authorcr Pawan Goyal}
\affil[$\ddagger\ddagger$]{ 
	\email{goyalp@mpi-magdeburg.mpg.de}, \orcid{0000-0003-3072-7780}}

\author[$\mathparagraph$]{Peter Benner}
\affil[$\mathparagraph$]{ 
	\email{benner@mpi-magdeburg.mpg.de}, \orcid{0000-0003-3362-4103}}

\shorttitle{GN-SINDy for Discovery of Partial Differential Equations}
\shortauthor{Forootani et al.}
\shortdate{}


\msc{MSC1, MSC2, MSC3}
  
\abstract{
The sparse identification of nonlinear dynamical systems (\sindy) is a data-driven technique employed for uncovering and representing the fundamental dynamics of intricate systems based on observational data. However, a primary obstacle in the discovery of models for nonlinear partial differential equations (\pde s) lies in addressing the challenges posed by the curse of dimensionality and large datasets. Consequently, the strategic selection of the most informative samples within a given dataset plays a crucial role in reducing computational costs and enhancing the effectiveness of \sindy-based algorithms. To this aim, we employ a greedy sampling approach to the snapshot matrix of a \pde~to obtain its valuable samples, which are suitable to train a deep neural network (\dnn) in a \sindy~framework. \sindy~based algorithms often consist of a data collection unit, constructing a dictionary of basis functions, computing the time derivative, and solving a sparse identification problem which ends to regularised least squares minimization. In this paper, we extend the results of a \sindy~based deep learning model discovery (\deepymod) approach by integrating greedy sampling technique in its data collection unit and new sparsity promoting algorithms in the least squares minimization unit. In this regard we introduce the greedy sampling neural network in sparse identification of nonlinear partial differential equations (\gnsindy) which blends a greedy sampling method, the \dnn, and the \sindy~algorithm. In the implementation phase, to show the effectiveness of \gnsindy, we compare its results with \deepymod~by using a \texttt{Python} package that is prepared for this purpose on numerous \pde~discovery. 
}

\novelty {
We introduce the Greedy Sampling Neural Network in Sparse Identification of Nonlinear Partial Differential Equations (\gnsindy), a pioneering approach that seamlessly integrates a novel greedy sampling technique, deep neural networks, and advanced sparsity-promoting algorithms. Our method not only addresses the formidable challenges posed by the curse of dimensionality and large datasets in discovering models for nonlinear \pde s but also sets a new standard for efficiency and accuracy by redefining the data collection and minimization units within the \sindy~framework. By combining the strengths of these diverse techniques, \gnsindy~represents a leap forward in the realm of model discovery, promising unprecedented insights into the intricate dynamics of complex systems.
}

\maketitle

  
\section{Introduction}%
\label{sec:intro}

Nonlinear dynamical systems are often encountered in various scientific fields, ranging from physics and biology to economics and engineering. Such systems can be extremely complex and difficult to understand, especially when they involve a large number of variables. Sparse identification of nonlinear dynamics (\sindy \cite{brunton2016sparse}) is a powerful approach that can help unravel the mysteries of these systems. Using a combination of machine learning and optimisation techniques, \sindy~can identify the governing equations of a nonlinear dynamical system from noisy and scarce data \cite{forootani2023, goyal2022neural}.

When dealing with partial differential equations (\pde s), \sindy~can be adapted to identify sparse representations of the nonlinear terms in the \pde~\cite{rudy2017data, both2021deepmod}. Traditional methods for \pde~identification often depend on a combination of domain specific knowledge, mathematical derivations, and experimental data. However, these approaches can be labor-intensive and may encounter limitations when applied to complex or poorly understood systems. Overcoming these challenges requires innovative techniques that can efficiently identify \pde s, provide a more comprehensive understanding of intricate physical phenomena, and enable improved predictions in diverse applications \cite{rudy2019data, tanyu2023deep}.


In the existing body of literature, the predominant focus has been on solving \pde s through either analytical or numerical means. Analytical approaches involve techniques like variable transformations to render the equation amenable or the derivation of an integral form of the solution \cite{de2006integral}. While these methods find applicability in handling straightforward \pde s, their efficacy diminishes when faced with more intricate equations. On the other hand, numerical methods aim to approximate the solution of a \pde~by discretizing its domain and solving a set of algebraic equations. Widely used numerical techniques include the finite difference method (\fdm) \cite{zeneli2021numerical} and the finite element method (\fem) \cite{bai2022local}.


Besides solving the \pde s with conventional methods, recent advances in machine learning techniques have proved their potential to address \pde~problems in scenarios with limited dataset. This implies having access solely to the \pde~problem data, rather than an extensive set of value pairs for the independent and dependent variables \cite{blechschmidt2021three}. In addition, modern machine learning software environments have provided automatic differentiation capabilities for functions realized by deep neural networks (\dnn) which is a mesh-free approach and can break the curse of dimensionality of the conventional methods~\cite{lu2021deepxde}. This approach was introduced in \cite{raissi2019physics}, where the term physics-informed neural networks(\pinn s) was coined. With the emergence of \pinn, employing a neural network has become a prominent method to construct a surrogate for data and subsequently conduct sparse regression on the network's predictions \cite{schaeffer2017learning, berg2019data, both2021model}. Alternatively, Neural \ode s have been introduced to unveil unknown governing equations \cite{rackauckas2020universal} from physical datasets. Diverse optimization strategies, employing the method of alternating direction, are explored in \cite{chen2021physics}, and graph-based approaches are formulated in \cite{seo2019differentiable}. Symmetry incorporation into neural networks is addressed directly by \cite{cranmer2020lagrangian}, utilizing both the Hamiltonian and Lagrangian frameworks. Furthermore, auto-encoders have been utilized to model \pde s and uncover latent variables \cite{iten2020discovering}. However, this approach does not yield an explicit equation and demands substantial amounts of data. It is worth to highlight that unlike the traditional \pde~solvers that focus more on methods such as the \fdm~\cite{zeneli2021numerical} and \fem~\cite{bai2022local}, the \dnn~based approaches (such as \pinn) are mesh free and therefore highly flexible.

The technique that we focus in this paper for \pde~identification is based on \sindy~algorithm, i.e. a method able to select, from a large dictionary, the correct linear, nonlinear, and spatial derivative terms, resulting in the identification associated \pde s from data \cite{brunton2016discovering}. In \sindy~only those dictionary terms that are most informative about the dynamics are selected as part of the discovered \pde. Previous sparse identification algorithms faced a number of challenges \cite{brunton2016discovering}: They were not able to handle sub-sampled spatial data, and the algorithm did not scale well to high-dimensional measurements. Standard model reduction techniques such as proper orthogonal decomposition (\propod) were used to overcome the high-dimensional measurements, allowing for a lower-order \ode~model to be constructed on energetic \propod~modes. This procedure resembles standard Galerkin projection onto \propod~modes \cite{rudy2017data}. \sindy~already applied on various model discovery applications, including for reduced-order models of fluid dynamics \cite{loiseau2018constrained} and plasma dynamics \cite{kaptanoglu2021physics}, turbulence closures \cite{beetham2020formulating}, mesoscale ocean closures \cite{zanna2020data}, nonlinear optics \cite{sorokina2016sparse}, computational chemistry \cite{boninsegna2018sparse}, and numerical integration schemes \cite{thaler2019sparse}. However, \sindy~algorithm in its original form has drawbacks such as sensitivity to accurate derivative information, lack of performance in the scarce data, sensitivity to the noisy dataset \cite{rudy2017data}.

Enhancements to the \sindy~framework have aimed at improvement its resilience to noise, providing uncertainty quantification, and adapting it for the modeling of stochastic dynamics \cite{boninsegna2018sparse, niven2020bayesian, messenger2021weak, hirsh2022sparsifying, wang2022data}. Nevertheless, these extensions have typically depended on computationally intensive methods for acquiring the necessary knowledge about probability distributions.


\subsection{Contribution}\label{contribution}
Expressing the unknown differential equation as $\partial_t \bu =\mathbf{f}(\mathbf{u}, \mathbf{u}_x, \dots)$ and assuming that the right-hand side is a linear combination of predefined terms, i.e., $\mathbf{f}(\mathbf{u},\mathbf{u}_x,\dots)= a \mathbf{u} + b \mathbf{u}_x + \dots = \Theta \xi$, simplifies model discovery to identifying a sparse coefficient vector $\xi$. The challenge lies in computing the time derivative $\bu_t$ and the function dictionary $\Theta$, especially when dealing with large dataset. The associated error in these terms tends to be high due to the utilization of numerical differentiation techniques like finite difference or spline interpolation. This limitation confines classical model discovery to datasets characterized by low noise and dense sampling. In contrast, deep learning-based methods overcome this challenge by constructing a surrogate from the data and determining the feature library $\Theta$, along with the time derivative $\bu_t$, through automatic differentiation. The \dnn~can be integrated seamlessly into \sindy~algorithm and be employed to effectively model the data, therefore facilitating the construction of a comprehensive function dictionary. 

A pivotal aspect of this approach involves the dynamic application of a mask during training loop of the \dnn, selectively activating terms in the function dictionary. In this regard the \dnn~is constrained to conform to solutions derived from these active terms. To determine this mask, any non-differentiable sparsity-promoting algorithm can be employed, such as \texttt{STRidge} \cite{brunton2016discovering}. This sophisticated approach not only allows for the use of a constrained neural network to precisely model the data and construct a robust function library but also employs an advanced sparsity-promoting algorithm to dynamically unveil the underlying equation based on the network's output. However, as the the number of training samples increase for a \dnn, the longer it takes to train the network. To alleviate this situation, it becomes crucial to choose a set of informative samples for training the network together with sparsity promoting algorithm (\sindy), and later, for \pde~model discovery.


In this study, the Discrete Empirical Interpolation Method (\deim~  \cite{drmac2016new}) is employed to strategically select informative samples from a snapshot matrix associated to a \pde. The objective is to reduce the dimensionality of the high-dimensional system while retaining essential features, ultimately enhancing the efficiency of subsequent data-driven methodologies. \deim~is a discrete variant of the empirical interpolation method (\eim) designed for constructing approximations of non-affine parameterized functions in continuous bounded domains, offering an associated error bound on the quality of approximation \cite{barrault2004empirical}. In particular, we utilize \qdeim~method which improves upon the original \deim~algorithm by providing a superior upper bound error and exhibiting numerically robust, high-performance procedures available in widely-used software packages such as \texttt{Python}, \texttt{LAPACK}, \texttt{ScaLAPACK}, and \texttt{MATLAB} \cite{drmac2016new}.

The selected \qdeim~samples are employed to train a \dnn~within the context of a \sindy-based algorithm. This integrated approach leverages the strengths of both \qdeim~and neural networks, where \qdeim~efficiently captures the most relevant information from the \pde~solution snapshots, and the neural network learns the underlying dynamics. Therefore, we name our proposed approach greedy sampling neural network for sparse identification of nonlinear dynamical systems (\gnsindy). In a detailed exploration, we systematically assess the influence of most informative samples obtained through \qdeim~on the snapshot matrix of a partial differential equation (\pde) for the purpose of discovering its governing equation. Additionally, we conduct a comparative analysis between the outcomes of \gnsindy~and \deepymod \cite{both2021deepmod} which uses random sampling approach. Our results underscore the substantial effect of leveraging highly informative samples on the training process of the \dnn~architecture for sparse identification of nonlinear \pde. Moreover, a \texttt{Python} package is prepared to support different implementation phases of \gnsindy~and its comparison with \deepymod~on \pde~model discovery corresponding to the Allen-Cahn equation, Burgers' equation and Korteweg-de Vries equation.


This article is organised as follows: In \Cref{sec: sindybackground}, we provide an overview of the \sindy~algorithm, focusing on its application for identifying \pde s. We elaborate on our approach to greedy sampling and its utilization in selecting the most informative samples from the dataset associated with \pde s in \Cref{sec:qdeim}. The \gnsindy algorithm, introduced in \Cref{sec:method}, is dedicated to the discovery of \pde s. In \Cref{sec:Simulation}, we delve into the simulation and comparative analysis of \gspinn~with \deepymod. The paper concludes with a summary in \Cref{sec:conclusion}.


\section{An Overview of \sindy~for \pde~identification}%
\label{sec: sindybackground}

Consider a nonlinear partial differential equations of the form
\begin{equation}\label{pde_1}
\mathbf{u}_t = \mathbf{f}(\mathbf{u},\mathbf{u}_x, \mathbf{u}_{xx}, \dots, x),
\end{equation}
where the subscripts denote partial differentiation in either time $t \in [0,t_{\maxm} ]$, or
space $x \in [x_{\minm}, x_{\maxm}] $, and $\mathbf{f}$ is an unknown right-hand side that is generally a non-linear function of $\mathbf{u}(x, t)$, and its derivatives. The goal in \sindy~algorithm is to discover $\mathbf{f}(\cdot)$ from the time series measurements of the dynamical system at a set of spatial points $x$. The main assumption behind \sindy~is that the un-known function $\mathbf{f}(\cdot)$ has only a few terms which makes the functional form sparse with respect to large space of possible contributing terms. For instance we can name Allen-Chan equation ($\mathbf{f}= 0.0001\mathbf{u}_{xx} - 5\mathbf{u}^3  + 5\mathbf{u}$) or Korteweg-De Vries equation ($ \mathbf{f} = -6 \mathbf{u} \mathbf{u}_x - \mathbf{u}_{xxx} $).

\sindy~algorithm is summarized in the following four steps.

\begin{enumerate}
    \item Data collection:

To do so we first measure $\mathbf{u}$ at $m$ different time instances and $n$ spatial locations, and we construct a single column vector $\bU \in \R^{n\cdot m}$.

\item Constructing the dictionary of basis functions:

At the second step we have to construct a dictionary consist of $D$ candidate liner, nonlinear, and partial derivatives
for the \pde:
\begin{equation}\label{dictionary_2}
\Theta(\bU)= \big[\mathbf{1},\ \bU,\ \bU^2,\ \dots,\ \bU_x,\ \bU \bU_x,\ \dots \big], \ \ \Theta \in \R^{nm \times D},
\end{equation}
where each single column of the matrix $\Theta$ contains all of the values of a particular candidate function across all of the $n\cdot m$ space-time grid points that data are collected. For example if we measure our process at $200$ spacial locations, at $300$ time instances and we include $40$ candidate terms in the \pde, then $\Theta \in \R^{200 \cdot 300 \times 40}$. These basis functions represent the possible dynamics that the system can exhibit. The choice of basis functions depends on the specific problem and can range from simple polynomials to more complex functions such as trigonometric functions or exponential functions and it is quite essential for the success of the algorithm. 

\item Computing the time derivative:

The third step in the \sindy~algorithm is to compute the time derivative of $\bU$, which is often implemented numerically. Having $\bU_t$ and other ingredients we can write the system of equation \eqref{pde_1} in the following form:
\begin{equation}\label{pde_4}
\bU_t = \Theta(\bU) \xi,
\end{equation}
where $\xi \in \R^{D}$ is an unknown vector that has to be computed by proper algorithm and its elements are coefficients corresponding to the active terms in the dictionary $\Theta(\cdot)$ that describe the evolution of the dynamic system in time. We show by $\xi_i$, the $i$th element of this column vector.

\item Solving the sparse identification problem:
The last and fourth step in the \sindy~algorithm is to compute the coefficient vector $\xi$ through a least squares optimization formulation which is a Nondeterministic Polynomial time (NP) mathematical hard problem. Therefore, there is a need to introduce a regularization technique that promotes sparsity. In the realm of regression analysis, different techniques address the challenges associated with model complexity and overfitting. Ordinary least squares (\texttt{OLS}) optimization, a classical method, often leads to complex models vulnerable to overfitting. Recognizing this limitation, sparse optimization techniques, such as \texttt{LASSO} and \texttt{Ridge} regression, introduce regularization to the \texttt{OLS} objective function \cite{friedman2001elements,tibshirani1996regression}. \texttt{LASSO} employs an $l_1$ penalty, promoting sparsity in the coefficient vector and facilitating effective feature selection, while \texttt{Ridge} regression, utilizing an $l_2$ penalty, enhances stability by addressing multicollinearity and preventing overfitting.

The sequential threshold ridge regression (\stridge) algorithm builds on these principles, introducing a regularized variant of \ols~that effectively deals with challenges in discovering physical laws within highly correlated and high-dimensional datasets. \stridge~addresses the limitations of standard regression methods, making it particularly useful in scenarios involving spatio-temporal data or complex, correlated features. Each of these approaches represents a nuanced trade-off between model complexity and predictive accuracy, allowing researchers and practitioners to tailor their choices based on the specific characteristics of their datasets and the goals of their modeling endeavors.

In the original \sindy~algorithm authors used \stridge~to unveil elusive governing equations, typically expressed as \ode s \cite{brunton2016discovering}. This approach laid the foundation for a broader framework known as \pde-functional identification of nonlinear dynamics (\pde-\texttt{FIND})\cite{rudy2017data}, facilitating the discovery of unknown relationships within the system \cite{vaddireddy2020feature}. In some other approaches as discussed in \cite{beck2013sparsity,yang2016sparse}, the authors considered to have some prior knowledge regarding non-zero terms in the coefficient vector $\xi$, which contradicts with the main assumption of \sindy~based approaches. The hyperparameters in \stridge~include the regularization weight $\lambda_{\texttt{STR}}$ tolerance level $\texttt{tol}_{\texttt{STR}}$, and maximum iteration $\texttt{max-iter}_{\texttt{STR}}$ which are
to be tuned to identify appropriate physical models. The convergence of analysis of sequential thresholding can be found in \cite{zhang2019convergence}. 

\end{enumerate}

\section{Importance of the sampling methods for \pde~discovery}\label{sec:qdeim}
By strategically opting the most informative samples, we can ensure that the identified \pde~faithfully represents the intricate dynamics of the system. This careful selection not only enhances the predictive capacity of the model but also provides valuable insights into the underlying physical processes. Moreover, choosing the most informative samples expedite the training time of the \dnn, and hence reduces the associated computational cost. Conversely, an inadequate sampling strategy may lead to inaccuracies in \pde~identification, hindering scientific advancements in diverse fields. Therefore, the judicious combination of an effective sampling strategy with the consideration of informative samples is pivotal for the successful and reliable identification of \pde s.

Within the existing body of literature, several residual point sampling methodologies have predominantly found application. Noteworthy among these are non-adaptive uniform sampling techniques, including (i) the equispaced uniform grid, (ii) uniformly random sampling, (iii) Latin hypercube sampling \cite{raissi2019physics}, (iv) the Halton sequence \cite{halton1960efficiency}, (v) Hammersley sequence~\cite{Wongetal97}, and (vi) the Sobol sequence \cite{sobol1967distribution}. Additionally, adaptive non-uniform sampling approaches, such as Residual-based Adaptive Distribution (\texttt{RAD}) and Residual-based Adaptive Refinement with Distribution (\texttt{RAR-D}) \cite{wu2023comprehensive}, have been explored. Despite the potential offered by these methodologies, they exhibit a strong dependence on specific problem characteristics and commonly entail laborious and time-intensive processes. Notably, adaptive sampling methods like \texttt{RAD} or \texttt{RAR-D} necessitate dataset re-sampling within the training loop of the \nn, introducing considerable computational overhead to the associated algorithm. For a comprehensive exploration of non-adaptive and residual-based adaptive sampling strategies and their comparisons in \pinn~training, the interested reader is directed to the detailed investigation presented in \cite{wu2023comprehensive}.


In numerous practical scenarios, like phase-field modeling or fluid dynamics, sensors are commonly employed to measure state variables. The quantity of sensors is often restricted by physical or financial constraints, and strategically situating these sensors is vital for attaining accurate estimations. Unfortunately, identifying the optimal locations for sensors to deduce the parameters of the \pde~poses an inherent combinatorial challenge. Existing approximation algorithms may not consistently yield efficient solutions across all pertinent cases. The matter of optimal sensor placement captivates research attention, extending its significance even into domains such as control theory and signal processing~\cite{both2021model}.

Literature have documented five categories of strategies for positioning sensors, as follows: (i) techniques relying on proper orthogonal decomposition (\propod) \cite{zhang2016pod} or sparse sensing like compressed sensing \cite{moslemi2023sparse}, (ii) methodologies involving convex optimization \cite{joshi2008sensor}, (iii) algorithms guided by a greedy approach, exemplified by Frame-Sens \cite{ranieri2014near}, (iv) heuristic methods, including population-based search \cite{lau2008tabu}, and (v) application of machine learning methodologies \cite{wang2019reinforcement}. The applicability of these methods is limited to particular cases due to conservative assumptions, complexity in implementation for large scale problems.


An interesting numerical technique that can be employed as a sampling method is \qdeim~algorithm which is originated in the context of model order reduction for high-dimensional dynamical systems \cite{chaturantabut2010nonlinear}. It is particularly useful when dealing with large-scale problems, such as those arising in computational physics or engineering simulations. \qdeim~aims to identify a reduced set of basis functions that capture the essential dynamics of the system by exploiting the empirical interpolation idea. The method tactically selects a sparse set of interpolation points from the high-dimensional state space, enabling an efficient representation of the system's behavior. \qdeim~has proven to be effective in reducing the computational cost associated with solving complex systems by constructing a low-dimensional surrogate model while preserving key system features. 

In the case of \pde~discovery, and specifically \dnn~framework, with an increase in the number of training samples, the training duration also extends. To address this issue, selecting a set of informative samples for network training becomes crucial. Our subsequent goal is to identify the governing \pde~underlying the data set within \sindy~procedure. To be more precise, we methodically investigate the impact of the most informative samples acquired via \qdeim~on \pde~snapshot matrix.

\subsection{Notes on \qdeim~algorithm}
Herein, we provide a concise introduction to \qdeim~algorithm, emphasizing its association with a well-known model order reduction technique known as \propod. \propod~has been used widely to select measurements in the state space that are informative for feature space reconstruction \cite{manohar2018data}. Consider the set of snapshots $\{u_1,\dots,u_m\}\in \R^n$ and an associated snapshot matrix $\cU=[u_1,\dots,u_m] \in \R^{n \times m}$ that is constructed by measuring the solution at $m$ different time instances and $n$ different spatial locations of a \pde. In the conventional \propod, we construct an orthogonal basis that can represent the dominant characteristics of the space of expected solutions that is defined as $\text{Range}\ \cU$, i.e., the span of the snapshots. We compute the singular value decomposition (\SVD) of the snapshot matrix $\cU$,
\begin{equation}\label{svd_1}
\cU = \bZ  \mathbf{\Sigma} \bY^\top,
\end{equation}
where $\bZ \in \R^{n\times n}$, $\mathbf{\Sigma} \in \R^{n \times m}$, and $\bY \in \R^{m \times m}$ with $\bZ^\top \bZ = \bI_n$, $\bY^\top \bY = \bI_m$, and $\mathbf{\Sigma}=\text{diag}(\sigma_1,\sigma_2,\cdots,\sigma_z)$ with $\sigma_1\ge \sigma_2 \ge \dots \ge \sigma_z \ge 0$ and $z = \min\{m,n\}$. The \propod~will select $\bV$ as the leading $\texttt{r}$ left singular vectors of $\bU$ corresponding to the $r$ largest singular values. Using \texttt{Python-Numpy} array notation, we denote this as $\bV = \bZ [:,:\texttt{r}]$. The basis selection via \propod~minimizes $ \bV:= \min_{\Phi \in \R^{n\times r}} \| \cU - \Phi \Phi^\top \cU \|^2_F$, where $\|\cdot \|_F$ is the Frobenius norm, over all $\Phi \in \R^{n \times r}$ with orthonormal columns. In this regard, we can say $\cU = \bZ  \mathbf{\Sigma} \bY^\top \approx \bZ_{\texttt{r}} \mathbf{\Sigma}_{\texttt{r}} \bY_{\texttt{r}}^\top $, where matrices $\bZ_{\texttt{r}}$ and $\bY_r^\top$ contain the first $\texttt{r}$ columns of $\bZ$ and $\bY^\top$, and $\mathbf{\Sigma}_{\texttt{r}}$ contains the first $\texttt{r} \times \texttt{r}$ block of $\mathbf{\Sigma}$. More details regarding \propod~can be found in \cite{kunisch2002galerkin}.


While the reduced-order model resides within the $\texttt{r}$-dimensional subspace, the conventional \propod~encounters a challenge when transitioning back to the original space. To address this issue, various approaches have been proposed in the literature, with \deim~being one such solution \cite{chaturantabut2010nonlinear}. Notably, \deim~offers the distinct advantage of flexibility, allowing its outcomes to extend beyond the realm of model order reduction, especially in the context of nonlinear function approximation. Additionally, the performance of the original \deim~algorithm has been enhanced by incorporating \texttt{QR}-factorization, resulting in two notable improvements: (i) a reduction in upper bound error and (ii) increased simplicity and robustness in implementation. \qdeim~leverages the pivoted factorization of \texttt{QR} factorization and the \SVD, providing a robust sampling method. Specifically, we employ \texttt{QR} factorization with column pivoting on $\bZ^\top_\texttt{r}$ and $\bY^\top_\texttt{r}$ to identify the most informative samples in the snapshot matrix $\cU$. This pivoting technique offers an approximate greedy solution for feature selection, termed submatrix volume maximization, as the matrix volume is defined by the absolute value of the determinant. Note that \texttt{QR}-factorization has been implemented and optimized in most scientific computing packages and libraries, such as \texttt{MATLAB}, and \texttt{Python}. Further details about \qdeim, its theoretical analysis, and applications can be found in \cite{manohar2018data, chaturantabut2010nonlinear}.

\subsection{Applying \qdeim~algorithm on \pde~dataset}
We employ the \qdeim~algorithm to acquire the most valuable samples on the spatio-temporal grid, utilizing the snapshot matrix $\cU$. To accomplish this, initially \texttt{SVD} is performed on the snapshot matrix $\cU$, yielding matrices $\bZ$, $\mathbf{\Sigma}$, and $\bY^\top$. Choosing $\texttt{r}$ leading singular values is determined based on a precision value $\epsilon_{\texttt{thr}}$, which is associated with the underlying dynamical system and is subject to heuristic selection. In the literature, $\epsilon_{\texttt{thr}}$ is commonly known as the energy criterion \cite{gavish2014optimal}. Specifically, the $\texttt{r}$ leading singular values are chosen to satisfy the following criterion:
\begin{equation}
1 - \frac{ \sum_{j=1}^{\texttt{r}} \sigma_j }{ \sum_{k=1}^{z} \sigma_k } < \epsilon_{\texttt{thr}} ,\ \  \texttt{r} < z,
\end{equation} 
once the desired level of precision is attained, we create a reduced approximation of the snapshot matrix $\cU$ by extracting the first $\texttt{r}$ columns of matrix $\bZ$, the first $\texttt{r}$ singular values from the diagonal matrix $\mathbf{\Sigma}$, and the initial $\texttt{r}$ rows of $\bY^\top$. In \texttt{Python} notation, this process is expressed as $\bZ_\texttt{r} = \bZ[:,:\texttt{r}]$, $\mathbf{\Sigma}\texttt{r} = \mathbf{\Sigma}[:\texttt{r},:\texttt{r}]$, and $\bY^\top\texttt{r} = \bY^\top[:\texttt{r},:]$. To identify significant time and space indices, we implement \texttt{QR} decomposition with column pivoting on the reduced-order matrices $\bY^\top_\texttt{r}$ and $\bZ^\top_\texttt{r}$, housing the foremost $\texttt{r}$ left and right singular vectors. The indices corresponding to the most informative spatio-temporal points in the snapshot matrix $\cU$ are denoted as $\texttt{ind}_x$ and $\texttt{ind}t$. For simplicity, we represent pairs of space-time points as $(t_i,x_i)$ and their associated solutions as $u_i$. \Cref{diemalgo} succinctly outlines the essential steps of the \qdeim~sampling approach, taking the snapshot matrix $\cU$, spatio-temporal domains $x,\ t$, and the precision value $\epsilon_{\texttt{thr}}$ as inputs and yielding sampled pairs $(t_i, x_i)$ along with the corresponding measured values $u(t_i,x_i)$. The cardinality of the sampled dataset is denoted as $\cN$.


\subsubsection{Importance of domain division}\label{time_division}

To analyze the local dynamics in the dataset and pinpoint optimal points in the spatio-temporal domain, we utilize a method that involves partitioning the time domain into uniform intervals. This approach employs the \qdeim~on each sub-domain to enhance efficiency. The parameter representing the number of divisions in the time domain is denoted as $\texttt{t}_{\texttt{div}}$. In this particular setup, the subdomains are non-overlapping, resulting in the total number of selected points being a combination of samples chosen from each subdomain.

In the context of \texttt{Python} notation, if we opt to divide the time domain into three parts, the representation will be as follows: the first part as $\cU[:,:\texttt{m}/3]$, the second part as $\cU[:,\texttt{m}/3:2\texttt{m}/3]$, and the third part as $\cU[:,2\texttt{m}/3:]$. This decomposition strategy is designed to capture the local dynamics of the partial differential equation (\pde) in each subdomain, focusing on sampling the most informative segments of the snapshot matrix. The rationale behind this lies in recognizing that the system's behavior might showcase variations across different domains, with specific physical attributes exhibiting notable distinctions. Noteworthy examples include issues related to abrupt features like shock waves, which may manifest these differences. Furthermore, dividing a large domain into smaller sub-domains and independently applying \qdeim~to each sub-domain serves to alleviate the necessity for intricate neural network structures in \pde~discovery. This streamlined approach enhances computational efficiency and facilitates a more targeted analysis of local variations within the dataset.

\begin{algorithm}[t]
	\KwData{$\cU,\ \{x_k\}_{k=1}^{n},\ \{t_k\}_{k=1}^{m} ,\ \epsilon_{\texttt{thr}}$.}
	\KwResult{ $\texttt{ind}_t$, $\texttt{ind}_x$, domain sampled pairs $(t_i, x_i)$ and $u(t_i,x_i)$.
	}
	$\texttt{r} = 1 $\;
	
	$\bZ, \mathbf{\Sigma}, \bY^\top \gets$ \texttt{SVD}($\cU$), \Comment computing \texttt{SVD} on snapshot matrix $\cU$\;
	
	Find the lowest $ \texttt{r}$, such that $ 1 - \frac{\sum_{j=1}^{\texttt{r}} \sigma_j}{\sum_{j=1}^{z} \sigma_j}  \geq \epsilon_{\texttt{thr}}$;
	
	$\bZ_\texttt{r} \gets \bZ[:, : \texttt{r}]$, $\bY^\top_\texttt{r} \gets \bY^\top[:\texttt{r}, :]$, \Comment{selecting $\texttt{r}$ dominant left and right singular vectors}; 
	
	$\texttt{ind}_x \gets$ \texttt{QR}($\bZ^\top_\texttt{r}$, pivoting = \texttt{True}), \Comment{storing pivots from pivoted \texttt{QR} factorization of $\bZ^\top_\texttt{r}$};
	
	$\texttt{ind}_t \gets$ \texttt{QR}($\bY_\texttt{r}$, pivoting = \texttt{True}), \Comment{ storing pivots from pivoted \texttt{QR} factorization of $\bY^\top_\texttt{r}$};
	
	$x_i \gets \text{from}\ \texttt{ind}_x $, $t_i \gets \text{from}\ \texttt{ind}_t$, $u(t_i,x_i)$;
	
	\caption{Sample selection based on a two-way \qdeim~procedure}
	\label{diemalgo}
\end{algorithm}

\section{\gnsindy: Greedy Sampling Neural Network for Sparse Identification of \pde s}%
\label{sec:method}

The \dnn~has been widely used in literature for solving \pde s due to their strength as the universal function approximators that can represent and learn complex relationships between input data and output solutions. This property is especially valuable in the context of \pde s, where finding analytical solutions may be challenging or impossible for certain complex systems \cite{huang2022partial}. However, for the case of \pde~model discovering from a given snapshot matrix there exist a few works that addressed such problem. Our \dnn~structure is similar to the work reported in \cite{both2021deepmod} with mainly difference in the \qdeim~algorithm that we employ to select the most valuable samples in the training loop. Several other works extended the results of \cite{both2021deepmod} for the case of noisy and scarce dataset by using integration scheme in the \dnn~training loop \cite{forootani2023, goyal2022neural}.


The process of model discovery through deep learning involves the utilization of a neural network to generate a surrogate model, denoted as $\hat{u}$, for the given data $u$. A collection of potential terms, represented by the dictionary $\Theta$, is established through automatic differentiation from $\hat{u}$. The neural network is then restricted to solutions permissible within this term dictionary. The network's loss function encompasses two key components: (i) a mean square error aimed at learning the mapping $(t,x)\to \hat{\bu}$ and (ii) a term designed to impose constraints on the network's solutions,
\begin{equation}\label{loss}
\mathbf{\mathcal{L}} = \frac{1}{\cN} \sum_{i=1}^{\cN} \Big( {\bu(t_i,x_i) - \hat{\bu}(t_i,x_i)}\Big)^2 + \frac{1}{\cN} \sum_{i=1}^{\cN} \Big( \frac{\partial \hat{\bu}(t_i,x_i)}{\partial t_i} - \mathbf{\Theta}\big( \hat{\bu}(t_i,x_i) \big) \xi \Big)^2,
\end{equation}
where the domain sampled pairs $(t_i,x_i)$, and $u(t_i,x_i)$ are computed by employing \qdeim~algorithm on the snapshot matrix $\cU$. 

The sparse vector $\xi$ is being learned during the training loop, and has two roles: (i) identifying the active elements, or those with non-zero values, in the \pde, and (ii) imposing constraints on the network based on these active terms. In our implementation we dissociate these two objectives by uncoupling the constraint from the actual process of selecting sparsity. Initially, we compute a sparsity mask, denoted as $\texttt{g}$, and subsequently restrict the network solely based on the active terms identified within this mask. In the sense that rather than imposing constraints on the neural network using $\xi$ alone, we opt to constrain it with the element-wise multiplication of $\xi$ and $\texttt{g}$, thereby replacing equation \eqref{loss} with:
\begin{equation}\label{loss_mask}
\mathbf{\mathcal{L}} = \frac{1}{\cN} \sum_{i=1}^{\cN} \Big( {\bu(t_i,x_i) - \hat{\bu}(t_i,x_i)}\Big)^2 + \frac{1}{\cN} \sum_{i=1}^{\cN} \Big( \frac{\partial \hat{\bu}(t_i,x_i)}{\partial t_i} - \mathbf{\Theta}\big( \hat{\bu}(t_i,x_i) \big) (\xi \cdot \texttt{g}) \Big)^2.
\end{equation}

The training process according to equation \eqref{loss_mask} involves a two-step procedure. Initially, we compute the sparsity mask, denoted as $\texttt{g}$, utilizing a sparse estimator. Subsequently, we minimize it with respect to network parameters, employing the masked coefficient vector. Importantly, the sparsity mask $\texttt{g}$ does not necessitate a differentiable calculation, allowing for the utilization of any traditional, non-differentiable sparse estimator.

This approach has the following advantages: i) it furnishes an impartial estimate of the coefficient vector by abstaining from applying $l_1$ or $l_2$ regularization on $\xi$ ii) the sparsity pattern is derived from the complete dictionary $\Theta$ rather than solely from the persisting active terms, enabling the dynamic inclusion and exclusion of active terms throughout the training process. In this regard our training loop consist of a (i)function approximator that creates a surrogate model of the dataset, (ii) a dictionary $\Theta$ of possible terms and time derivatives $\frac{\partial\hat{u}}{\partial t}$, (iii) a sparsity estimator that creates a mask to choose the active columns in the dictionary using sparse regression technique, and (iv) a constraint that imposes the function approximator to the solutions allowed by the active terms obtained by sparsity estimator.

\paragraph{Notes on sparsity mask and sparsity estimator}
The sparsity mask $\texttt{g}$ is computed using an estimator which is not involved in the training loop of the \dnn~structure and it is not differentiable. Therefore, the right procedure to update the sparsity mask $\texttt{g}$ within the training loop is essential for the success of the algorithm. In particular, it is important to allow the function approximator unit learns the solution $\bu$ for some iterations, then the sparsity mask be updated. We control the update of the sparsity mask $\texttt{g}$ with parameters \texttt{patience}, \texttt{periodicity}. We define the \texttt{patience} as the number of iterations that we allow the \dnn~to learn the solution $\bu$ before applying the sparsity mask. Moreover, the $\texttt{periodicity}$ refers to the consistent intervals at which we check whether sparsity mask needs to be updated or not within the training iterations. Sparsity threshold $\delta_{\texttt{spr}}$ is a value that we impose on the solution that is computed via sparsity estimator to update the sparsity mask $\texttt{g}$. Note that the sparse estimator solves $ \frac{\partial \hat{\bu}(t_i,x_i)}{\partial t_i} - \mathbf{\Theta}\big( \hat{\bu}(t_i,x_i) \big) \xi_{\texttt{est}} = 0 $ , where $ \xi_{\texttt{est}} $ is its solution. The elements of the vector $ \xi_{\texttt{est}}$ that have absolute values less than sparsity threshold $\delta_{\texttt{spr}}$ will be used to update the sparsity mask $\texttt{g}$. It is worth to highlight that an appropriate algorithm that promotes sparsity can be used to compute $\xi_{\texttt{est}}$ such as \stridge~\cite{rudy2017data}, \lasso~\cite{friedman2010regularization}.

\begin{figure}[t]
	\centering
	\includegraphics[width=\textwidth]{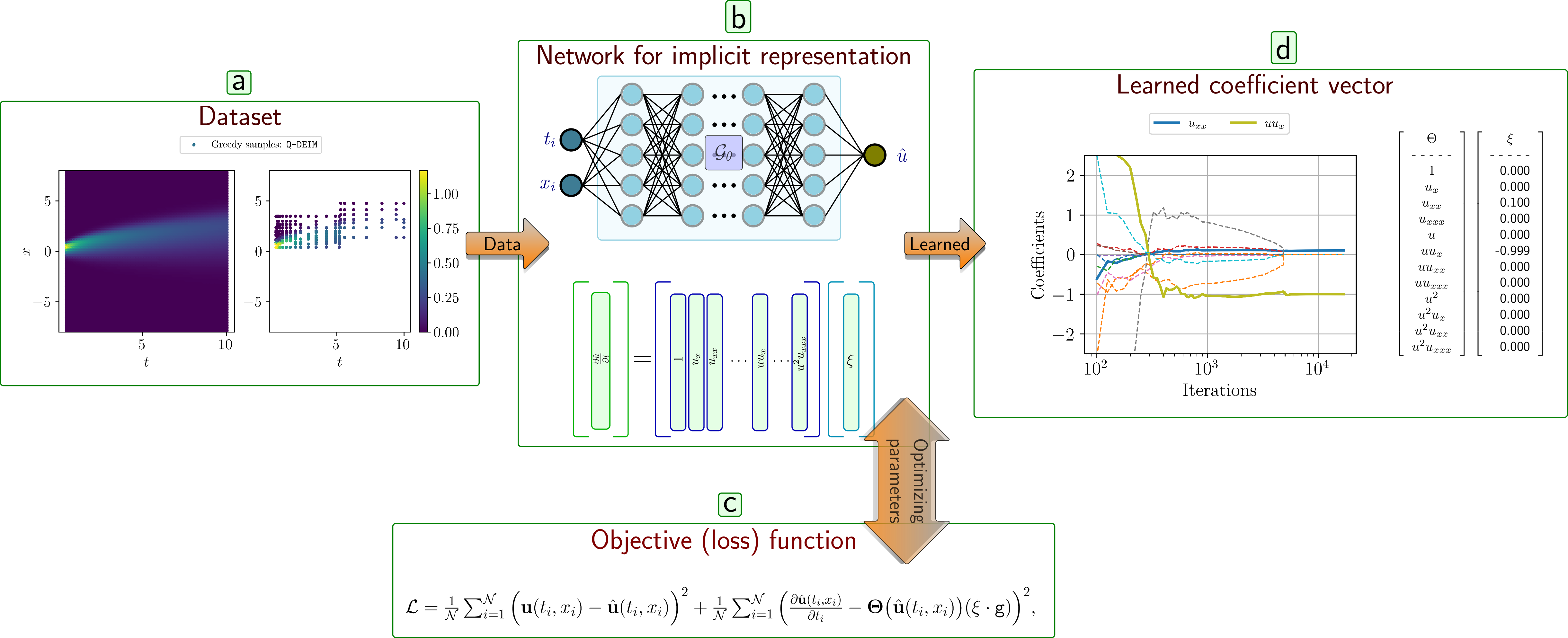}
	\caption{
    A schematic diagram of \gnsindy~for \pde~discovery. (a) sampling \pde~dataset with \qdeim~algorithm to choose the most informative samples, (b) feeding the
    valuable sample pairs $(t_i,x_i)$ computed by \qdeim~into \dnn~and utilizing the output of the \dnn~as the function approximator $\hat{u}$ to construct the dictionary $\Theta$ and to compute the $\frac{\partial \hat{u}}{\partial t}$ via automatic differentiation (c) estimating the coefficient vector $\xi$ via stochastic gradient descent subject to the loss function $\cL$ and sparsity mask.}
 \label{gsindy_schematic}	
\end{figure}

\paragraph{Notes on constraint}
The sparse coefficient vector $\xi$ is computed concurrently through the \dnn~training loop. When the sparsity mask $\texttt{g}$ is updated, the dictionary columns corresponding to inactive terms are omitted in the solution of $ \frac{\partial \hat{\bu}( \cdot )}{\partial t} - \mathbf{\Theta}\big( \hat{\bu}(\cdot) \big) \xi =0 $. Note that by default it is assumed that all the dictionary terms are active (sparsity mask $\texttt{g}$). Omitting inactive terms in the dictionary $\Theta$ in parallel to decreasing mean square error which decreases the complexity of finding the coefficient vector $\xi$. Moreover, any suitable algorithm that promotes sparsity can be used to compute $\xi$ such as \stridge~\cite{rudy2017data}, \lasso~\cite{friedman2010regularization}.\Cref{alg_two} summarizes the procedure that has been explained in this section. Moreover a schematic diagram of \gnsindy~algorithm is depicted in \Cref{gsindy_schematic}.


\alglanguage{pseudocode}
\begin{algorithm}[h]
\KwData{$\cU$, $x$, $t$,  $\epsilon_{\texttt{thr}}$ for the \qdeim~ algorithm, a neural network $\cG_\theta$ (parameterized by $\theta$), maximum iterations ${\texttt{max-iter}}$, \texttt{patience}, \texttt{periodicity}, and sparsity threshold $\delta_{\texttt{spr}}$}

\KwResult{ Estimated coefficient vector $\xi$}

$(t_i, x_i),\ u(t_i,x_i) \gets$ Apply \qdeim($\cU$) based on \Cref{diemalgo} \Comment{selecting most informative samples}\;

Initialize the \dnn~module parameters\; 
$k=1$\;
Initialize the sparsity mask $\texttt{g}$\;

\While{ $k<$ \texttt{max-iter} }{
- Feed the domain pairs $(t_i,x_i)$ as an input to the \dnn~($\cG_\theta$) and predict output $\mathbf{\hat{u}}(t_i,x_i)$\;

- Construct the library $\Theta \big(\mathbf{\hat{u}}(t_i,x_i) \big)$\;

- Compute the derivative information $\frac{\partial \mathbf{\hat{u}}(t_i,x_i)}{ \partial t_i }$ using automatic differentiation\;

- Compute the coefficient vector $\xi$ by sparsity promoting algorithm such as \texttt{STRidge}, subject to the sparsity mask $\texttt{g}$ to solve $ \frac{\partial \hat{\bu}(t_i,x_i)}{\partial t_i} - \mathbf{\Theta}\big( \hat{\bu}(t_i,x_i) \big) (\xi \cdot \texttt{g} )=0 $\;

- Compute the cost function~\eqref{loss_mask}\;

- Update the parameters of $\cG_\theta$ and the coefficient vector $\xi$ using gradient descent\;

\If{ (k - \texttt{patience}) \%\texttt{periodicity} == 0  }{
        - Update the sparsity mask $\texttt{g}$ using \texttt{Lasso} or \texttt{STRidge} algorithm subject to the solution $ \frac{\partial \hat{\bu}(t_i,x_i)}{\partial t_i} - \mathbf{\Theta}\big( \hat{\bu}(t_i,x_i) \big) \xi_{\texttt{est}} =0 $ and sparsity threshold $\delta_{\texttt{spr} }$\;  \Comment{note that updating the sparsity mask is done independent from estimated coefficients}
    }

}

\caption{GN-SINDy: greedy sampling neural network for sparse identification of nonlinear \pde s }
\label{alg_two}
\end{algorithm}
	
\section{Simulation Results} \label{sec:Simulation}
In this section we employ \gnsindy~algorithm on several \pde s including, Burgers' equation, Allen-Cahn equation, and Korteweg-de Vries equation. We also compare the results of the \gnsindy~with \texttt{DeePyMoD} \cite{both2021deepmod} where samples are selected randomly for \pde~ discovery. The simulation examples have different level of complexity and non-linearity. We emphasize that, within our framework, we do not take into account the boundary conditions of the \pde s. This is because our primary objective is not to solve these equations but rather to discover the governing equation. Each example undergoes a two-part simulation process. In the first part, we apply \qdeim~to the provided snapshot matrix of a \pde~to pinpoint the most informative samples.
In the second part, we make use of most informative samples to train our \dnn~ and estimate the coefficient vector $\xi$.

Since we have no prior knowledge of the various hyperparameters required for \gnsindy, we do some trial and error to find the right configuration setup. The \qdeim~ algorithm, which filters out informative samples from unimportant ones, is the main step in our procedure. A suitable trade-off must be considered for the time division $\texttt{t}_{\texttt{div}}$ and the precision value $\epsilon_{\texttt{thr}}$, since the cardinality of the sample set can be varied by setting these parameters. We first fix the time division parameter $\texttt{t}_{\texttt{div}}$ for each example and vary the precision value $\epsilon_{\texttt{thr}}$ to obtain the most informative samples corresponding to each example. Then, we randomly select a part of these informative samples to train the \dnn~structure for \pde~descovery. It is worth highlighting that the cardinality of the candidate informative samples used in the \gnsindy~algorithm is considered as the minimum possible value needed to recover the \pde. We identify this minimum number by testing different setups for \qdeim~algorithm, sparse estimator, constraint, sparsity scheduler and \dnn~structure. We show that using the \qdeim~algorithm on the \pde~snapshot matrix $\cU$ provides us with valuable insight to capture the local dynamics, so that even randomly selecting the fraction of resulting subsamples will increase the success rate of discovering the governing equation compared to the approaches that directly randomly sample the \pde~snapshot matrix, such as \deepymod.

\paragraph{Data generation.}
The dataset in this article are taken from repositories reported in \cite{rudy2017data, both2021deepmod}. Our data preprocessing step only includes selecting the most informative samples by employing \qdeim~algorithm on the snapshot matrix of \pde~dataset. Our procedure to select the samples has two stages: (i) we employ \qdeim~with the appropriate setting for the precision threshold $\epsilon_{\texttt{thr}}$ and the time division $\texttt{t}_{\texttt{div}}$, (ii) we choose randomly portion of the \qdeim~samples for the training loop of the \gnsindy. For each \pde~we also mention the original range of the space-time domain dataset.

\paragraph{Architecture.}
Our method utilizes multi-layer perceptron networks equipped with periodic activation functions, particularly adopting the \texttt{SIREN} architecture introduced by \cite{sitzmann2019siren}. This approach facilitates the extraction of implicit representations from measurement data, allowing us to customize the number of hidden layers and neurons for each specific example.
We specify number of hidden layers and associated number of neurons for each \pde~inclusively.

\paragraph{Hardware.}
Our endeavor to uncover governing equations via neural network training and parameter estimation necessitated the utilization of cutting-edge computational resources. To this end, we employed an \nvidia RTX A4000 GPU, with its robust $16$ GB of RAM, to handle the computationally demanding aspects of our research. For CPU-intensive tasks, such as data generation, we harnessed the power of a 12th Gen \intel~\coreifive-12600K processor, equipped with a remarkable $32$ GB of memory.
\paragraph{Training set-up.}
We use the Adam optimizer which is a popular optimization algorithm in deep learning. The optimizer is configured with $\texttt{learning\_rate}= 10^{-3}$, signifying the step size for updating the model's parameters during training \cite{paszke2017automatic}. Notably, the $\texttt{beta}$ values for the exponential moving averages were set to $0.99$, placing increased emphasis on past gradients in the optimization process \cite{paszke2017automatic}. Additionally, the \texttt{amsgrad} variant was enabled, ensuring the stability of the optimizer's update rule for the moving average of squared gradients \cite{paszke2017automatic}. These hyperparameter choices were made to enhance the convergence and performance of the neural network model, demonstrating a thoughtful consideration of optimization strategies in the pursuit of improved computational results. Specification of the dictionary $\Theta$, the sparsity estimator to update the sparsity mask $\texttt{g}$, and the constraint on the coefficient vector $\xi$ for the \pde~discovery together with their related parameters will be given inclusively. Furthermore, number of time domain division ($\texttt{t}_{\texttt{div}}$) to apply \qdeim~for each \pde~dataset, precision value ($\epsilon_{\texttt{thr}}$), and maximum number of iterations (${\texttt{max-iter}}$) will be mentioned separately for each example. Lastly, to ensure reproducibility in our experiments, random number generators were seeded using specific values. The \texttt{NumPy} library was initialized with the seed $42$, and the \texttt{PyTorch} library with the seed $50$ in all the simulation examples. This deliberate seeding allows for the precise replication of random processes, facilitating the verification and validation of our results by other researchers



\subsection{Burgers' equation}
\label{subsec:Burger}


Burgers' equation, named after the Dutch mathematician and physicist Jan Martinus Burgers, is a fundamental \pde~that arises in the study of fluid dynamics and nonlinear waves. Introduced in 1948, this equation represents a simplified model for one-dimensional, inviscid fluid flow with small amplitude and shallow water conditions. Burgers' equation combines elements of both the linear advection equation and the nonlinear conservation law, making it a versatile tool in the analysis of various physical phenomena. Burgers' equation is often expressed as:
\begin{equation}
\frac{\partial u}{\partial t} + \gamma u \frac{\partial u}{\partial x} + \nu \frac{\partial^2 u}{\partial x^2} = 0,
\end{equation}
where $u(x,t)$ represents the velocity field of the fluid at a spatial location $x$ and time $t$ and $\gamma$ and $\nu$ are constant with nominal values $\gamma=-1$ and $\nu=0.1$, respectively.  

The physical interpretation of Burgers' equation can be understood by considering the propagation of a wave in a fluid. As the wave propagates, the fluid particles experience shearing forces due to the presence of the wave. These shearing forces cause a change in the velocity of the fluid particles, which is governed by Burgers' equation.

It is worth to highlight that Burgers' equation has been applied to model a wide range of physical phenomena, including: fluid flows, traffic flow, and population dynamics \cite{kachroo2018traffic}. Moreover, it exhibits several important properties that makes it a versatile tool for modeling physical phenomena. In particular, the nonlinear term $u\frac{\partial u}{\partial x}$ gives rise to the formation of shock waves, which are characterized by sharp discontinuities in the velocity field. The presence of the diffusion term $\frac{\partial^2 u}{\partial x^2}$ imparts a dissipative nature to the equation. This means that the amplitude of the waves will decay over time, reflecting the energy dissipation due to viscosity. The equation also exhibits dispersion, meaning that the speed of the waves depends on their wavelength. This is a consequence of the linear advection term $u\frac{\partial u}{\partial x}$. These properties make Burgers' equation a powerful tool for analyzing a wide range of physical phenomena, from fluid flows to traffic patterns to population dynamics.

\paragraph{Analysis}
To generate the dataset for our simulation we use the solution of Burgers' equation\footnote{\url{ https://www.iist.ac.in/sites/default/files/people/IN08026/Burgers_equation_viscous.pdf}}. The spatial domain $x\in [-8,\ 8]$ and the time domain $t\in [0.5, 10]$ are considered each having $100$ samples, therefore the associate snapshot matrix $\cU \in \mathbb{R}^{100 \times 100}$. For the simulation setup of \gnsindy~we consider a $4$ layer \dnn~structure each layer having $64$ neurons, a dictionary consist of the combination of polynomial and derivative terms up to order $2$, a sparsity scheduler with $\texttt{periodicity}=100,\ \texttt{patience}=500$, a sparse estimator $\texttt{STRidge}$ with the hyperparameters maximum iteration $100$ and sparsity threshold $\delta_{\texttt{spr}}= 0.05$, a constraint $\texttt{STRidge}$ with maximum iteration $\texttt{max-iter}_{\texttt{STR}} = 100$ and tolerance $\texttt{tol}_{\texttt{STR}}=0.05$. For both sparse estimator and constraint the regularization weight set to $\lambda_{\texttt{STR}}=0$. Moreover, the maximum iteration of the training loop is considered $\texttt{max-iter}=25000$. With the mentioned settings our dictionary has the following terms
$$[1,\ \bu,\ \bu_{xx},\ \bu,\ \bu \bu_x, \bu \bu_{xx}, \bu^2, \bu^2 \bu_x, \bu^2 \bu_{xx}].$$

\begin{figure}[h]
\centering
\includegraphics[width=0.7\textwidth]{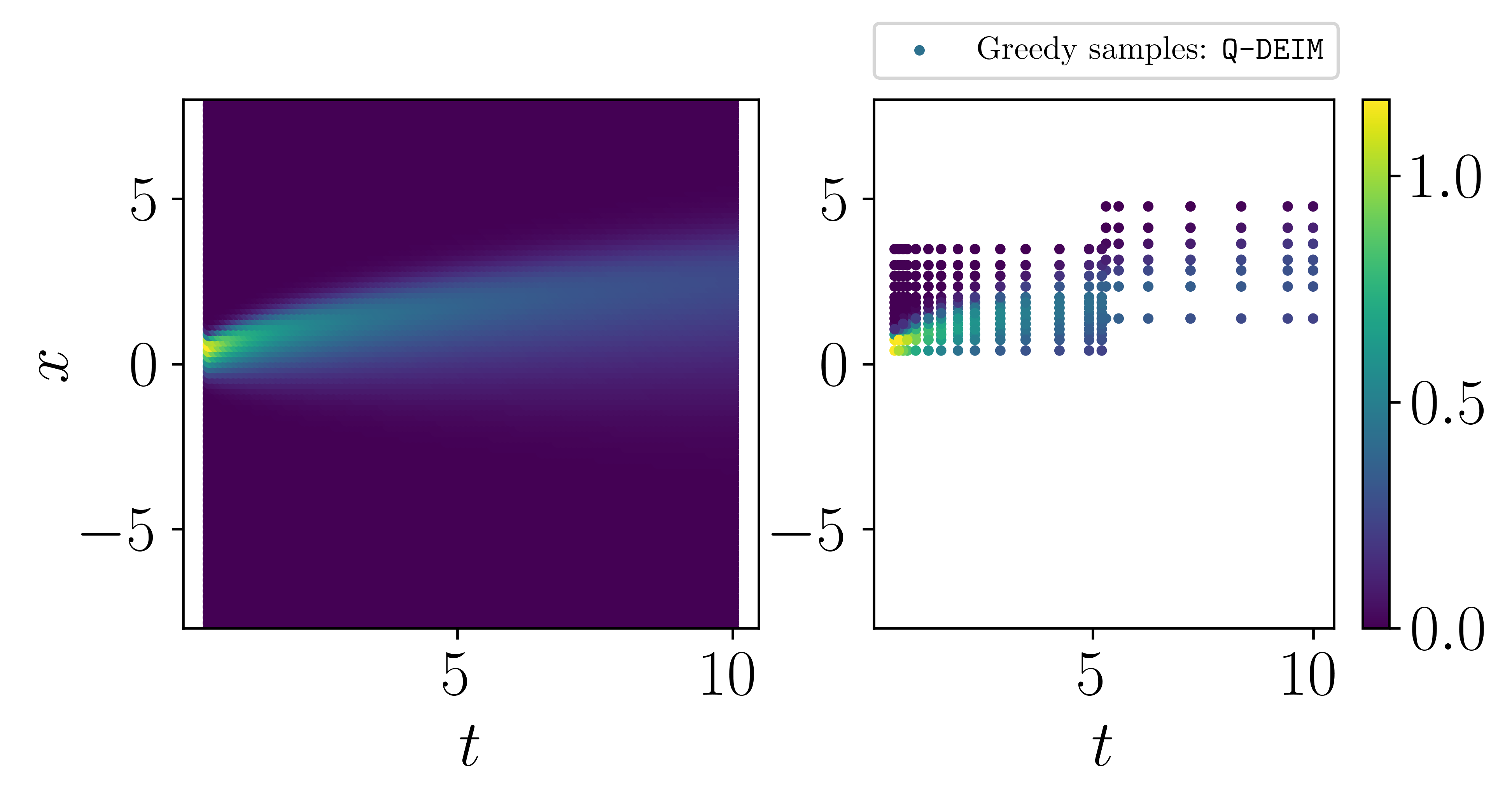}
\caption{(left) Entire dataset ; (right) Greedy samples resulted by \qdeim~algorithm for Burgers' equation with $\texttt{t}_{\texttt{div}}=2$ and $\epsilon_{\texttt{thr}}=10^{-5}$.}
\label{fig:burger_entire_deim}
\end{figure}

In a set of experiments, we consider to choose different values for the precision $\epsilon_{\texttt{thr}}$ of the \qdeim~algorithm with time division $\texttt{t}_{\texttt{div}}=2$ to see how it affects \gnsindy~performance. Set of different values are $\epsilon_{\texttt{thr}}=\{ 10^{-3},\ 10^{-4},\ 10^{-5},\ 10^{-6} \}$ which results $121,\ 180,\ 245$ and $1313$ samples respectively. With a few times trial and error we discover the Burgers' equation with $50$ random samples among total of $245$ samples corresponding to $\epsilon_{\texttt{thr}}= 10^{-5}$. Applying the \qdeim~algorithm with the mentioned setup on the snapshot matrix reveals that the valuable parts of the Burgers' equation lies on the interval $x\in [0,\ 5]$. In \Cref{fig:burger_entire_deim} we show the result of \qdeim~algorithm corresponding to $\texttt{t}_{\texttt{div}}=2$ and $\epsilon_{\texttt{thr}}=10^{-5}$.

To evaluate the \gnsindy~performance corresponding to each precision value $\epsilon_{\texttt{thr}}$ we choose $50$ random samples among the samples that \qdeim~ provided, for the training of our \dnn~structure. The results of the simulations is shown in the \Cref{table_Burger_precision_sensitivity}. This result proves the importance of choosing right precision value $\epsilon_{\texttt{thr}}$ to select the most informative samples. Decreasing precision value $\epsilon_{\texttt{thr}}$ forces \qdeim~algorithm to choose more samples form the snapshot matrix that these samples may not carry valuable information, therefore increases the failure rate of the \gnsindy. In addition, with the results of \Cref{table_Burger_precision_sensitivity}, we identify the proper setting for the \qdeim~algorithm hyperparameters which are $\epsilon_{\texttt{thr}}=10^{-5}$ and $\texttt{t}_{\texttt{div}}=2$. 

\begin{table}[h!]
  \begin{center}
    \caption{\gnsindy~performance with different precision value $\epsilon_{\texttt{thr}}$, fixed $\texttt{t}_{\texttt{div}}=2$ and fixed sample size $50$ for recovering Burgers' equation}
    \label{table_Burger_precision_sensitivity}
    \begin{tabular}{|c|c|c|}
      \toprule 
      \textbf{Precision value} & \textbf{Estimated \pde} \\
      \hline
      \midrule 

      $\epsilon_{\texttt{thr}}=10^{-3}$ &  $\bu_t + 0.0993 \bu_{xx}- 1.0037 \bu \bu_{xx} = 0$ \\ 
      $\epsilon_{\texttt{thr}}=10^{-4}$ &  $\bu_t + 0.1019 \bu_{xx}- 1.0236 \bu \bu_{xx} = 0$ \\ 

      $\epsilon_{\texttt{thr}}=10^{-5}$ & $\bu_t + 0.0985 \bu_{xx}- 0.9857 \bu \bu_{xx} = 0$ \\
      $\epsilon_{\texttt{thr}}=10^{-6}$ & $\bu_t + 0.0983 \bu_{xx}- 0.9969 \bu \bu_{xx} = 0$  \\
      \bottomrule 
    \end{tabular}
  \end{center}
\end{table}

In a more advanced analysis, we contemplate exploring how the performance of the \gnsindy~algorithm is influenced by varying sparse estimators and constraints. To do so, we maintain the simulation setup as previously specified and evaluate diverse combinations of sparse estimators and constraints, including \stridge, \lasso, and \ols. The sparsity threshold for sparse estimator of both types \stridge~and \lasso~is set to  $\delta_{\texttt{thr}}=0.05$. From the results reported in \Cref{Burger_sensitivity_estimator_constraint} we can conclude that choosing \stridge~as the sparsity estimator has a significant impact on recovering right \pde~model for Burgers' equation. Moreover, these results show that choosing most informative samples has to be combined with the suitable choice of sparse estimator and constraint so that \gnsindy~can recover the model correctly.

\begin{table}[h!]
  \begin{center}
    \caption{\gnsindy~performance with different sparse estimator and different constraint under fixed \qdeim~setting $\epsilon_{\texttt{thr}}=10^{-5}$, $\texttt{t}_{\texttt{div}}=2$ and fixed sample size $50$ for recovering Burgers' equation}
    \label{Burger_sensitivity_estimator_constraint}
    \begin{tabular}{|c|c|c|}
      \toprule 
      \textbf{Sparse estimator} & \textbf{Constraint} & \textbf{Estimated \pde}\\
      \hline
      \midrule 
      \stridge & \stridge & $\bu_t + 0.0985 \bu_{xx}- 0.9857 \bu \bu_{xx} = 0$\\
      \lasso & \stridge & $\bu_t + 0.1014 \bu_{xx} -0.7531 \bu \bu_x - 0.2754 \bu^2 \bu_{x} =0$\\
      \stridge & \ols & $\bu_t + 0.1014 \bu_{xx}- 1.0018 \bu \bu_{x}=0$\\
      \lasso & \ols & $\bu_t + 0.1011 \bu_{xx} -0.7481 \bu \bu_x - 0.2798 \bu^2 \bu_{x} =0$ \\
      \bottomrule 
    \end{tabular}
  \end{center}
\end{table}

In the next phase of simulation, we consider to evaluate the performance of \gnsindy~ under variation of the \dnn~structure. In this regard we fixed the previous setup, i.e. precision value $\epsilon_{\texttt{thr}}=10^{-5}$, $\texttt{t}_{\texttt{div}}=2$, selecting $50$ samples from samples earned by \qdeim~algorithm, sparse estimator and constraint both of the type \stridge. Moreover, the number of hidden layers is fixed to $4$ and we vary the number of neurons in each layer based on a geometric sequence with a common ratio of $2$ with initial value $8$ neurons. As reported in \Cref{table_Burger_NN_sensitivity}, the results show that all the configurations have a good performance and the combination of \stridge as the sparse estimator and as the constraint is robust under alteration of \dnn~structure except for the case of $8$ neurons in each layer.

\begin{table}[h!]
  \begin{center}
    \caption{\gnsindy~performance with different \dnn~structure, sparse estimator \stridge, constraint \stridge, fixed precision value $\epsilon_{\texttt{thr}}=10^{-5}$, fixed $\texttt{t}_{\texttt{div}}=2$ and fixed sample size $50$ for recovering Burgers' equation}
    \label{table_Burger_NN_sensitivity}
    \begin{tabular}{|c|c|c|}
      \toprule 
      \textbf{\dnn~structure} & \textbf{Estimated \pde} \\
      \hline
      \midrule 
      $[2,\ 8,\ 8,\ 8,\ 8,\ 1]$ & $\bu_t - 0.1014 \bu_{xx} -0.9054 \bu  \bu_x -0.2230 \bu^2 \bu_{xx}=0 $\\ 
      $[2,\ 16,\ 16,\ 16,\ 16,\ 1]$ & $\bu_t + 0.1006 \bu_{xx}- 1.0086 \bu \bu_{xx} = 0$  \\
      $[2,\ 32,\ 32,\ 32,\ 32,\ 1]$ &  $\bu_t + 0.1004 \bu_{xx}- 1.0091 \bu \bu_{xx} = 0$
      \\
       $[2,\ 64,\ 64,\ 64,\ 64,\ 1]$ & $\bu_t + 0.0985 \bu_{xx}- 0.9857 \bu \bu_{xx} = 0$ \\ 
      \bottomrule 
    \end{tabular}
  \end{center}
\end{table}



\begin{figure}[h!]
\centering
\begin{subfigure}{1\textwidth}
    \includegraphics[width=1\textwidth]{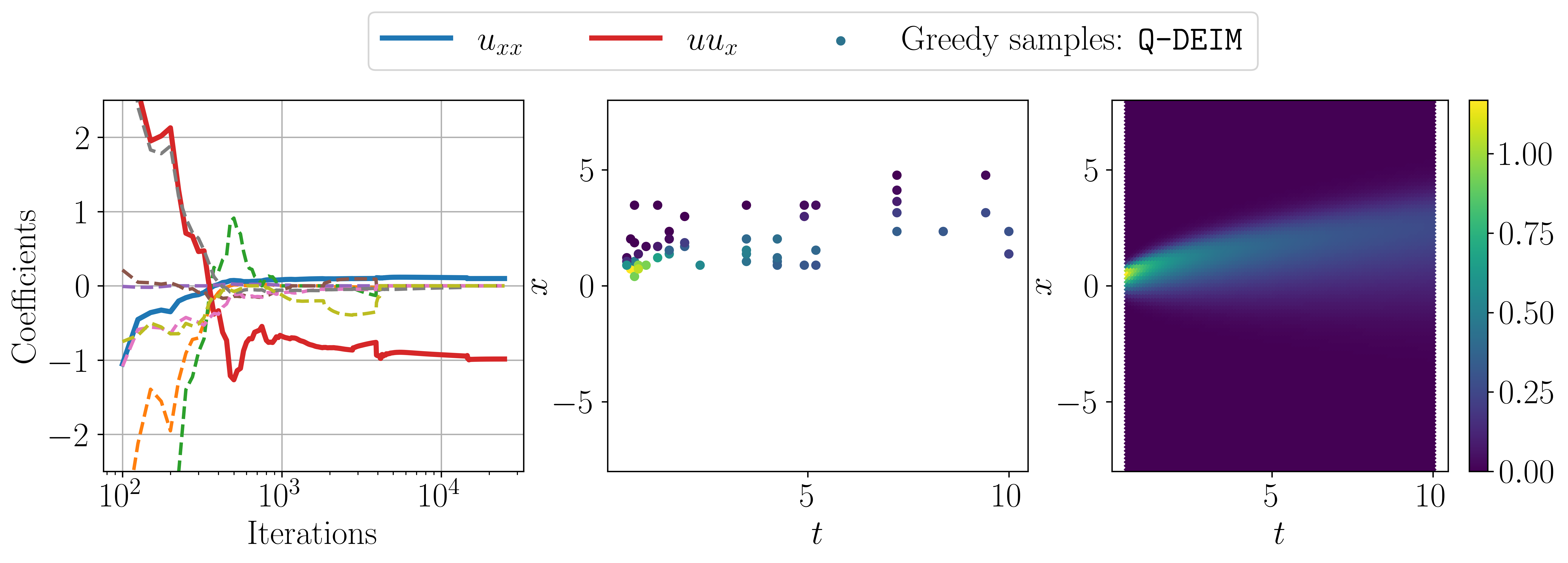}
    \caption{\gnsindy~performance: (right) Entire dataset; (middle) selection of $50$ Greedy samples resulted by \qdeim~algorithm for Burgers' equation with $\texttt{t}_{\texttt{div}}=2$ and $\epsilon_{\texttt{thr}}=10^{-5}$; (left) estimated coefficients with \gnsindy~through the iterations.}
    \label{fig:burger_deim_50}
\end{subfigure}
\hfill
\begin{subfigure}{1\textwidth}
    \includegraphics[width=1\textwidth]{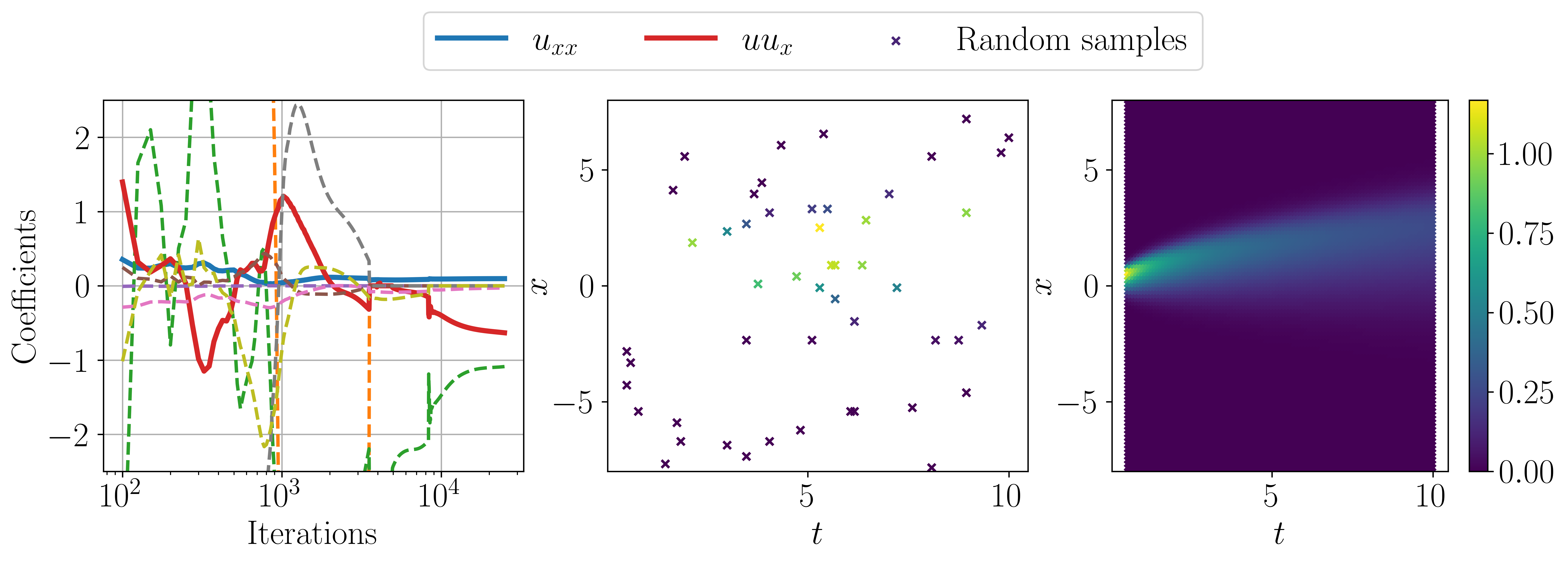}
    \caption{\deepymod~performance: (right) Entire dataset; (middle) selection of $50$ random samples for Burgers' equation; (left) estimated coefficients with \deepymod~through the iterations.}
    \label{fig:burger_random_50}
\end{subfigure}
\caption{Comparison of the \gnsindy~performance with \deepymod~in Burgers' equation model discovery}
\end{figure}

Finally, we compare the performance of \gnsindy~and \deepymod (\cite{both2021deepmod}). The same setup is assumed for both algorithms, with the knowledge that \deepymod~uses random sampling to choose samples from snapshot matrix, sparse estimator of the type \lasso~and constraint of the type \ols. In \Cref{fig:burger_deim_50}, and \Cref{fig:burger_random_50} the results of comparison between \gnsindy~ and \deepymod~are shown. It is quite obvious how \qdeim~algorithm assist \gnsindy~with $0.5\% (50/10000)$ of the entire dataset to successfully recover the right model. In particular, from \Cref{fig:burger_deim_50} (left) as the evolution of the coefficient demonstrates the $\bu_{xx}$ term is identified earlier while $\bu \bu_x$ has more fluctuations. Precise value of the estimated coefficients for Burgers' equation with sampling size $50$ and different methods are reported in \Cref{Burger_gnsindy_deepymod}. From this result we see that \gnsindy~outperforms \deepymod~and can recover the Burgers' equation with acceptable precision. 


\begin{table}[h!]
  \begin{center}
    \caption{Comparing \gnsindy~and \deepymod~for Burgers' \pde~discovery}
    \label{Burger_gnsindy_deepymod}
    \begin{tabular}{|c|c|c|}
      \toprule 
      \textbf{Algorithm} & \textbf{Estimated \pde} \\
      \hline
      \midrule 
     \gnsindy & $\bu_t + 0.0985 \bu_{xx} - 0.9857\bu \bu_x=0 $ \\
      \deepymod \cite{both2021deepmod} & $\bu_t -0.0287 \bu_x + 0.0971 \bu_{xx} -0.6326\bu \bu_x -1.0872\bu^2 \bu_x=0 $\\
      \bottomrule 
    \end{tabular}
  \end{center}
\end{table}

\subsection{Allen-Cahn equation}%
\label{subsec:AC}

The Allen-Cahn equation stands as a cornerstone of phase transition modeling, capturing the transformation of physical quantities, commonly referred to as order parameters, as a material undergoes transitions from one phase to another. Its elegance encompasses phenomena such as solidification, crystal growth, and the emergence of domain patterns in magnetic materials. The Allen-Cahn equation is typically expressed in one dimension as
\begin{equation}\label{AC_eq}
    \frac{\partial u }{\partial t} + \gamma_1 u_{xx} + \gamma_2 (u-u^3) = 0,
\end{equation}
where $u$ represents the order parameter, quantifying the level of order in the system, with nominal values $\gamma_1=0.0001$ and $\gamma_2 = 5$. For example in a binary alloy, $u=1$ corresponds to a pure $\texttt{A-phase}$, $u=-1$ corresponds to a pure $\texttt{B-phase}$, and $0$ represents a mixture of $\texttt{A}$ and $\texttt{B}$. The Allen-Cahn equation extends its reach beyond phase transitions, finding applications in various domains, such as magnetic Domain Formation \cite{shen2021cell}, phase Separation in Binary Alloys \cite{cahn1994evolution}, wetting of Surfaces \cite{zhang2023effect}, and pattern Formation in Biological Systems\cite{zhao2020learning}.

\paragraph{Analysis}

For Allen-Cahn equation we use the dataset that is reported in \cite{both2021deepmod}. The \pde~is discretized in $512$ spatial points and $201$ time instances, therefore the snapshot matrix $\cU \in \mathbb{R}^{512 \times 201}$ which demonstrates the curse of dimensionality regarding \pde~discovery. The setup of \gnsindy~algorithm for Allen-Cahn equation is as follows: a $4$ layer \dnn~structure each having $64$ neurons, a sparsity scheduler with $\texttt{patience}=1000$ and $\texttt{periodicity}=100$, both sparsity estimator and the constraint are of the type \stridge~with $\texttt{max-iter}_{\texttt{STR}}=100$ and tolerance $\texttt{tol}_{\texttt{STR}}=0.1$, a dictionary with polynomial and derivative terms up to order $3$. The maximum iteration for training loop is set to $\texttt{max-iter}=25000$. We note that the dictionary terms are 
$$[1,\  \bu_{x} ,\ \bu_{xx},\ \bu_{xxx},\ \bu ,\ \bu \bu_x,\ \bu \bu_{xx},\ \bu \bu_{xxx},\ \bu^2,\ \bu^2 \bu_x,\  \bu^2 \bu_{xx},\ \bu^2 \bu_{xxx},\ \bu^3,\ \bu^3 \bu_{x},\ \bu^3 \bu_{xx},\ \bu^3 \bu_{xxx} ].$$

\begin{figure}[tb]
\centering
\includegraphics[width=0.7\textwidth]{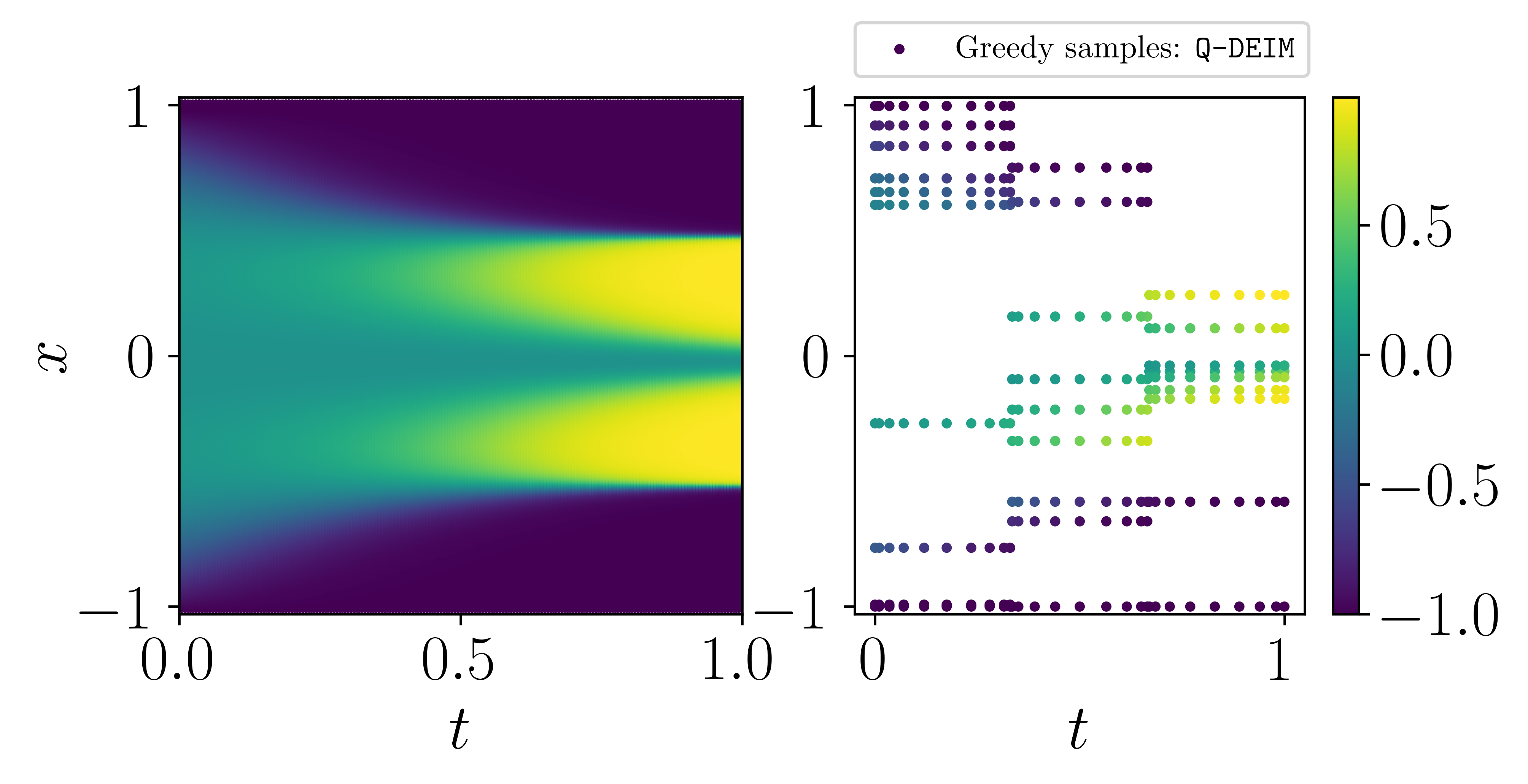}
\caption{(left) Entire dataset ; (right) Greedy samples resulted by \qdeim~algorithm for Allen-Cahn equation with $\texttt{t}_{\texttt{div}}=3$ and $\epsilon_{\texttt{thr}}=10^{-7}$.}
\label{fig:AC_entire_deim}
\end{figure}

To find out the minimum number of samples that is required to discover the Allen-Cahn equation we investigate the impact of precision value $\epsilon_{\texttt{thr}}$ on the performance of \gnsindy. In this regard we consider to evaluate \gnsindy~ with the following set of precision values $\epsilon_{\texttt{thr}}=\{10^{-5},\ 10^{-6},\ 10^{-7},\ 10^{-8}\}$. \qdeim~algorithm based on these precision values will result $147,\ 209$, $262$, and $386$ informative samples, respectively. With a few times trial and error we discover Allen-Cahn \pde~with 
$120$ samples out of total $262$ samples corresponding to the precision value $\epsilon_{\texttt{thr}} = 10^{-7}$. With this setting, these $262$ valuable samples are shown in \Cref{fig:AC_entire_deim}. From this figure we see how \qdeim~ algorithm selects informative samples to capture the local dynamics at each part of the domain. Moreover, it shows the importance of the time domain division for the sample selection of the \qdeim~algorithm.

To see how the other \gnsindy~configurations perform under different precision values we choose $120$ samples randomly out of $147,\ 209$, and $386$ and we evaluate the performance of the algorithm again. The results are reported in \Cref{AC_precision_sensitivity} which demonstrates that $\epsilon_{\texttt{thr}}=\ 10^{-7}$ has a better performance with respect to the other cases. Our purpose to do comparison for different precision value $\epsilon_{\texttt{thr}}$ is to find the best setting among all the possible combinations for the \pde~discovery. 


\begin{table}[h!]
  \begin{center}
    \caption{\gnsindy~performance with different precision value $\epsilon_{\texttt{thr}}$, $\texttt{t}_{\texttt{div}}=3$ and fixed sample size $120$ for recovering Allen-Cahn equation}
    \label{AC_precision_sensitivity}
    \begin{tabular}{|c|c|c|}
      \toprule 
      \textbf{Precision value} & \textbf{Estimated \pde} \\
      \hline
      \midrule 
      $\epsilon_{\texttt{thr}}=10^{-5}$ & $\bu_t + 0.0000 \bu_{xx}+ 4.8759 \bu -4.8517 \bu^3=0$ \\
      $\epsilon_{\texttt{thr}}=10^{-6}$ & $\bu_t + 0.0000 \bu_{xx}+ 4.9379 \bu -4.9026 \bu^3=0$ \\
      $\epsilon_{\texttt{thr}}=10^{-7}$ & $\bu_t + 0.0000 \bu_{xx}+ 4.9643 \bu -4.9328 \bu^3=0$ \\
      $\epsilon_{\texttt{thr}}=10^{-8}$ & $\bu_t + 0.0000 \bu_{xx}+ 4.8733 \bu -4.8171 \bu^3=0$ \\
      \bottomrule 
    \end{tabular}
  \end{center}
\end{table}


In the more quantitative analysis, we evaluate the performance of \gnsindy~under different choices of the sparse estimator and the constraint. For the \qdeim~settings we consider $\epsilon_{\texttt{thr}}=10^{-7}$, $\texttt{t}_{\texttt{div}}=3$ and we randomly select $120$ samples for the \dnn~training. In addition, the \dnn~structure is considered fixed as mentioned earlier, i.e. a $4$ layer \dnn with $64$ neurons in each layer. In \cref{AC_sensitivity_estimator_constraint} we report the results of these simulations. Almost all the configurations have the same performance, however the combination of \lasso~as the sparse estimator and \ols~as the constraint is slightly better than the others while the combination of \lasso~and \stridge~has .

\begin{table}[h!]
  \begin{center}
    \caption{\gnsindy~performance with different sparse estimator and different constraint under fixed \qdeim~setting $\epsilon_{\texttt{thr}}=10^{-7}$, $\texttt{t}_{\texttt{div}}=3$ and fixed sample size $120$ for recovering Allen-Cahn equation}
    \label{AC_sensitivity_estimator_constraint}
    \begin{tabular}{|c|c|c|}
      \toprule 
      \textbf{Sparse estimator} & \textbf{Constraint} & \textbf{Estimated \pde}\\
      \hline
      \midrule 
      \stridge & \stridge & $\bu_t + 0.0000 \bu_{xx}+ 4.9643 \bu -4.9328 \bu^3=0$\\
      \lasso & \stridge & $\bu_t + 0.0000 \bu_{xx}+ 4.8477 \bu -4.7882 \bu^3=0$\\
      \stridge & \ols & $\bu_t + 0.0000 \bu_{xx}+ 4.9790 \bu -4.9709 \bu^3=0$\\
      \lasso & \ols & $\bu_t + 0.0000 \bu_{xx}+ 5.0175 \bu -4.9724 \bu^3=0$\\
      \bottomrule 
    \end{tabular}
  \end{center}
\end{table}


We also consider to evaluate the \gnsindy~algorithm when we vary the structure of the \dnn. To do so, we fix the number of hidden layers to $4$ and we vary the number of neurons in each layer based on a geometric sequence with a common ration of $2$ with initial value $8$ neurons. The results of this experiments are reported in \Cref{table_AC_NN_sensitivity} where we \gnsindy is robust to the \dnn alteration.

\begin{table}[h!]
  \begin{center}
    \caption{\gnsindy~performance with different \dnn~structure, sparse estimator \stridge, constraint \stridge, fixed precision value $\epsilon_{\texttt{thr}}=10^{-7}$, fixed $\texttt{t}_{\texttt{div}}=3$ and fixed sample size $120$ for recovering Allen-Cahn equation}
    \label{table_AC_NN_sensitivity}
    \begin{tabular}{|c|c|c|}
      \toprule 
      \textbf{\dnn~structure} & \textbf{Estimated \pde} \\
      \hline
      \midrule 
      $[2,\ 8,\ 8,\ 8,\ 8,\ 1]$ &  $\bu_t + 0.0000 \bu_{xx}+ 5.0275 \bu -4.9655 \bu^3=0 $ \\
      $[2,\ 16,\ 16,\ 16,\ 16,\ 1]$ & $\bu_t + 0.0000 \bu_{xx}+ 4.9901 \bu -4.9428 \bu^3=0 $
      \\
      $[2,\ 32,\ 32,\ 32,\ 32,\ 1]$ &  $\bu_t + 0.0000 \bu_{xx}+ 4.9925 \bu -4.9481 \bu^3=0 $
      \\
       $[2,\ 64,\ 64,\ 64,\ 64,\ 1]$ & $\bu_t + 0.0000 \bu_{xx}+ 4.9643 \bu -4.9328 \bu^3=0 $ \\ 
      \bottomrule 
    \end{tabular}
  \end{center}
\end{table}

\begin{figure}[tb]
\centering
\begin{subfigure}{1\textwidth}
    \includegraphics[width=1\textwidth]{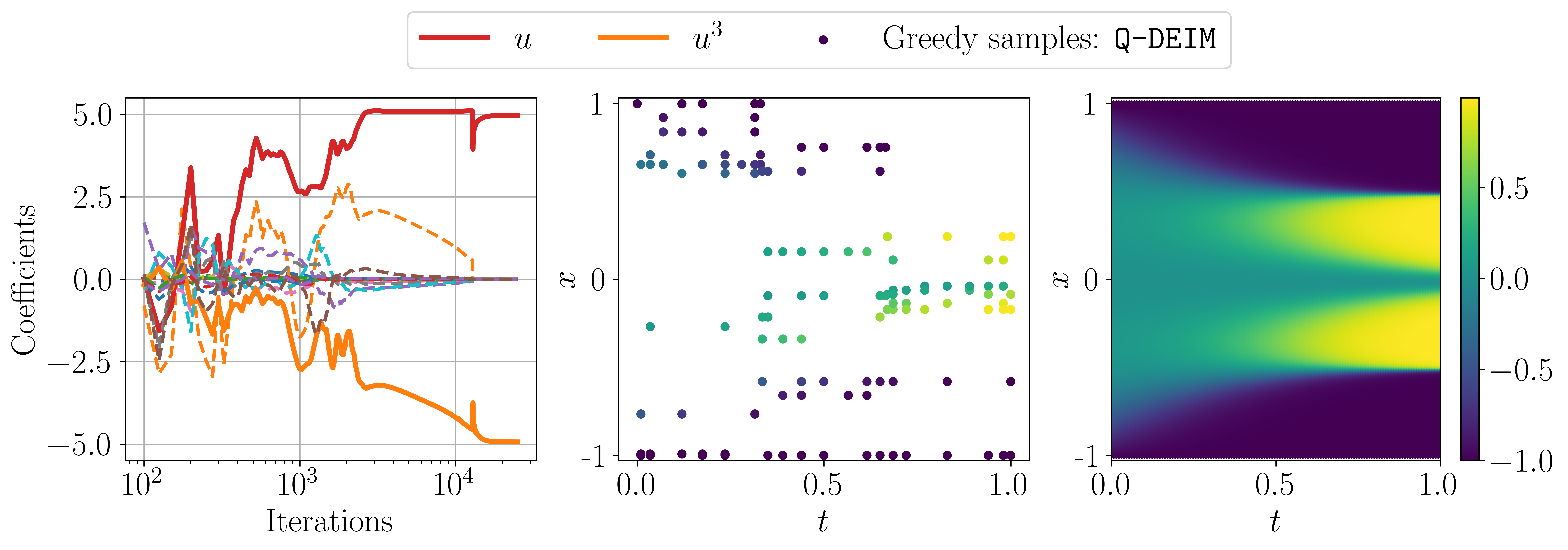}
    \caption{(right) Entire dataset ; (middle) selection of $120$ Greedy samples resulted by \qdeim~algorithm for Allen-Cahn equation with $\texttt{t}_{\texttt{div}}=3$ and $\epsilon_{\texttt{thr}}=10^{-7}$; (left) estimated coefficients with \gnsindy.}
    \label{fig:AC_deim_120}
\end{subfigure}
\hfill
\begin{subfigure}{1\textwidth}
    \includegraphics[width=1\textwidth]{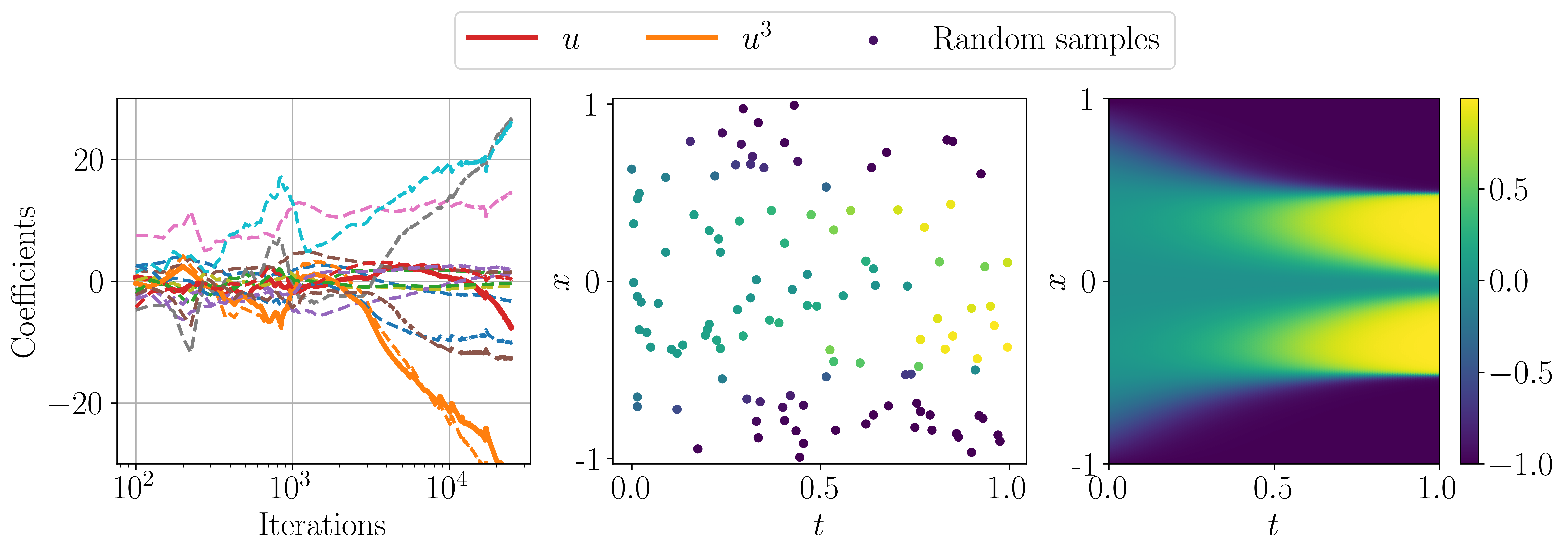}
    \caption{ (right) Entire dataset ; (middle) selection of $120$ random samples for Allen-Cahn equation; (left) estimated coefficients with \deepymod.
    }
    \label{fig:AC_random_120}
\end{subfigure}
\caption{ Comparison of the \gnsindy~performance with \deepymod~in Allen-Cahn equation model discovery}
\end{figure}

Moreover, we consider to do a comparison between \gnsindy~and \deepymod~(\cite{both2021deepmod}) to see how they 
perform in discovering the \pde~model with the same sample size. The same setup is considered for \deepymod~with its default sparse estimator as \lasso~and constraint of the type \ols. The results of the simulations are presented in \Cref{AC_gnsindy_deepymod}, that we see clearly \gnsindy~outperforms \deepymod~to recover the Allen-Cahn \pde~ with acceptable precision. It is worth to highlight that the coefficient corresponding to the term $\bu_{xx}$ is considerably less than the coefficient corresponding to the terms $\bu$ and $\bu^3$ in the ground truth equation, hence recovering such a small coefficient it is quite difficult for the \dnn. Moreover, in \Cref{fig:AC_deim_120} and \Cref{fig:AC_random_120} the selected greedy samples for \gnsindy, random samples for \deepymod~and evolution of different coefficients through the iteration are shown. The impact of \qdeim~algorithm in the success of recovering correct coefficients is quite clear. Vice versa \deepymod~has a weak performance due to the usage of randomly selected samples in the training loop of the \dnn~structure. 
\begin{table}[h!]
  \begin{center}
    \caption{\gnsindy~performance with different precision value $\epsilon_{\texttt{thr}}$, $\texttt{t}_{\texttt{div}}=3$ and fixed sample size $120$ for recovering Allen-Cahn equation}
    \label{AC_gnsindy_deepymod}
    \begin{tabular}{|c|c|c|}
      \toprule 
      \textbf{Algorithm} & \textbf{Estimated \pde} \\
      \hline
      \midrule 
     \gnsindy & $\bu_t + 0.0000 \bu_{xx}+ 4.9643 \bu -4.9328 \bu^3=0 $ \\
      \deepymod \cite{both2021deepmod}  & $\bu_t + 0.8913 + 14.6374 \bu_{x} - 3.2769 \bu_{xx} + 1.5190\bu_{xxx} - 7.7121 \bu 
 + 26.0523 \bu \bu_x$\\
 &$- 10.0496 \bu  \bu_{xx} + 1.4176 \bu \bu_{xxx} 
 - 43.0762 \bu + 26.6733 \bu^2 + 1.9071 \bu^2 \bu_x $\\ &$- 0.3155 \bu^2 \bu_{xx} - 32.1551 \bu^3 - 12.7150 \bu^3  \bu_{x} + 0.3304 \bu^3  \bu_{xx} - 0.8281 \bu^3  \bu_{xxx} = 0 $\\
      \bottomrule 
    \end{tabular}
  \end{center}
\end{table}


\subsection{Korteweg-de Vries (KdV) equation}\label{subsec:KdV}

The Korteweg-de Vries (KdV) equation is a nonlinear, dispersive partial differential equation that plays a central role in mathematical physics and nonlinear science. It was first introduced in 1895 by Dutch physicists Diederik Korteweg and Gustav de Vries to describe the propagation of shallow water waves. Since then, it has found applications in a wide range of physical systems, including plasma physics, optics, and condensed matter physics. The KdV equation is remarkable for its integrability, meaning that it can be solved exactly using various mathematical techniques. This property has led to the discovery of solitons, which are stable, localized waves that propagate without changing shape. Solitons have been observed in various physical systems, including shallow water waves, optical fibers, and Bose-Einstein condensates. The KdV equation has also been used to develop a deeper understanding of nonlinear dynamics. It has been shown that the KdV equation can exhibit a variety of complex behaviors, such as chaos and self-organization. These behaviors have relevance to a wide range of phenomena in physics, biology, and social sciences. The KdV equation remains a cornerstone of mathematical physics and nonlinear science, providing a rich and insightful framework for understanding nonlinear waves and their diverse applications \cite{yang2016exact}. It can be written as a partial differential equation in the form:
\begin{equation}\label{KDV_equation}
\frac{\partial u}{\partial t } + c\ u u_x + \alpha\ u_{xxx}=0, 
\end{equation}
where $c$ and $\alpha$ are constant variables with the nominal values $-6$ and $-1$, respectively.


\paragraph{Analysis}

The discretization of the KdV equation is outlined in \cite{rudy2017data}, which involves $512$ spatial points and $201$ temporal points within the original spatial range of $x \in [-30,\ 30]$ and a time span of $t \in [0,\ 20]$. The snapshot matrix $\cU \in \mathbb{R}^{512 \times 201}$ illustrates the challenge of dimensionality in training the \dnn~for \pde~discovery. In the first set of experiment to identify the most promising configuration of the \gnsindy~ to recover the KdV equation, we fix the time division $\texttt{t}_{\texttt{div}}=2$, and the following set of precision values are considered $\epsilon_{\texttt{thr}}=\{5 \times 10^{-5},\ 10^{-5},\ 10^{-6}, \  10^{-7}\}$. Now we need to set the \gnsindy~hyperparameters which include the following settings: to choose the polynomial order and the derivative orders for dictionary, to opt the type of sparse estimator and the constraint, and finally the setting corresponding to the \dnn~structure. A $4$ layer \dnn~structure is considered where each layer have $32$ neurons, sparse estimator and constraint of the type \stridge~with $\texttt{max-iter}_{\texttt{STR}}=100$ and tolerance $\texttt{tol}_{\texttt{STR}}=0.1$, sparse scheduler with $\texttt{patience}= 1000$ and $\texttt{periodicity}=50$, a dictionary with polynomial terms up to order $2$ and derivative terms up to order $3$. The maximum iteration of the \dnn~training loop is considered $\texttt{max-iter}=25000$. The \qdeim~algorithm provides $5725,\ 14450,\ 18432$ and $19801$ samples corresponding to each precision value. This set of filtered samples shows the difficulty of choosing proper number of informative samples for \pde~discovery. With the proposed settings for the order of polynomial as well as derivative terms of the dictionary the associated columns are   
$$ 
[1,\  \bu_{x} ,\ \bu_{xx},\ \bu_{xxx},\ \bu ,\ \bu \bu_x,\ \bu \bu_{xx},\ \bu \bu_{xxx},\ \bu^2,\ \bu^2 \bu_x,\  \bu^2 \bu_{xx},\ \bu^2 \bu_{xxx} ].
$$

\begin{figure}[tb]
\centering
\includegraphics[width=0.7\textwidth]{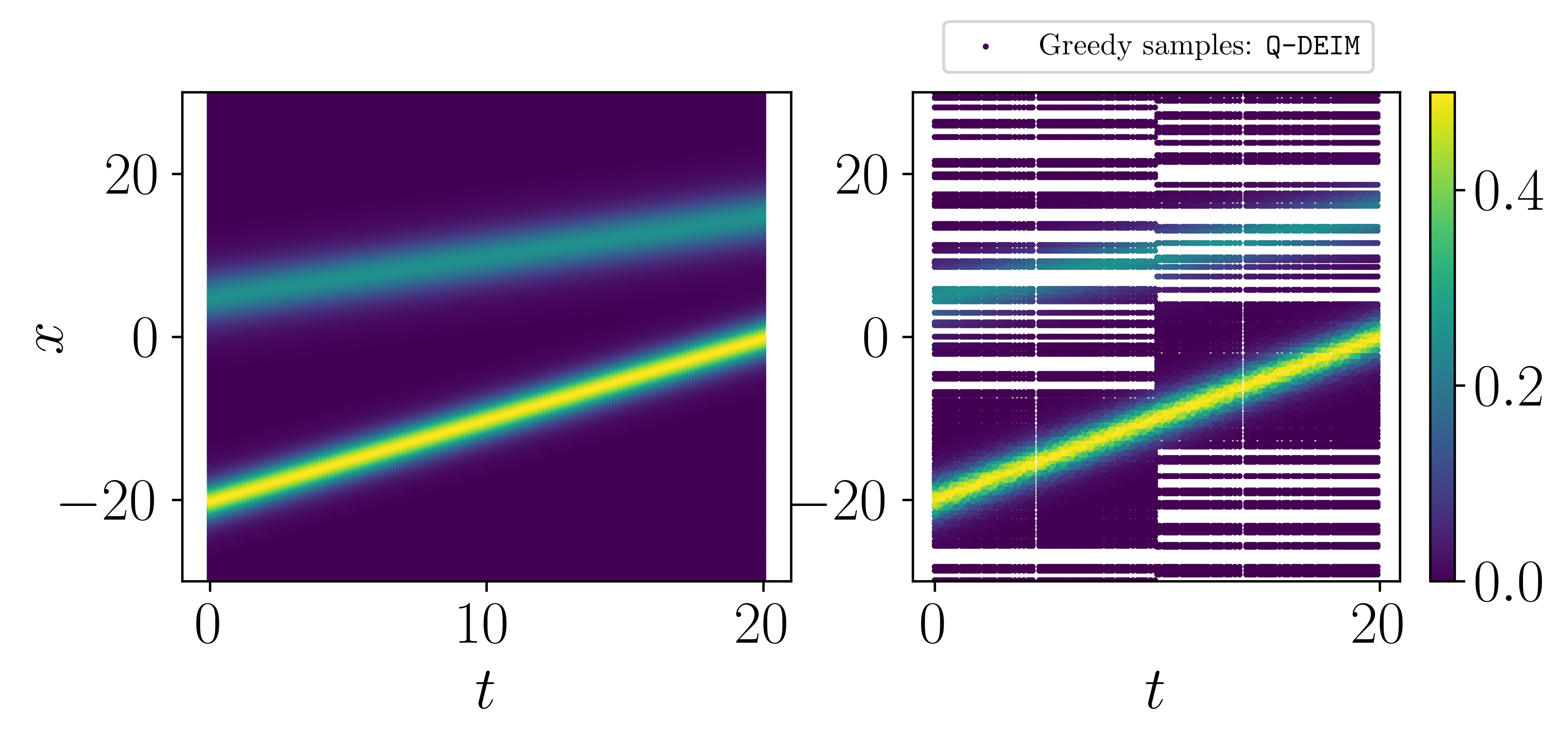}
\caption{(left) Entire dataset ; (right) Greedy samples resulted by \qdeim~algorithm for KdV equation with $\texttt{t}_{\texttt{div}}=2$ and $\epsilon_{\texttt{thr}}=10^{-5}$.}
\label{fig:KdV_entire_deim}
\end{figure}

With the mentioned \gnsindy~settings and among these set of informative samples corresponding to each precision value we randomly select portion of each dataset to identify the smallest possible number that can successfully recover the KdV equation. After a few times trial and error we realize to take $900$ samples from total number of informative samples corresponding to $\epsilon_{\texttt{thr}}=10^{-5}$ that results correct KdV equation. Note that we only use $0.87\%$ of the entire dataset approximately which shows significant reduction in computational cost. The outcomes of this experiment is shown in the \Cref{KdV_precision_sensitivity}, which demonstrate that choosing the precision value $\epsilon_{\texttt{thr}}=10^{-5}$ has the better performance with respect to the other choices. Moreover, these results prove the importance of choosing right precision value $\epsilon_{\texttt{thr}}$ in the success rate of \gnsindy~algorithm. Indeed decreasing $\epsilon_{\texttt{thr}}$ selects the samples that may not contribute significantly into the model discovery. Moreover, \Cref{fig:KdV_entire_deim} shows the entire dataset and selected samples by \qdeim~algorithm corresponding to $\texttt{t}_{\texttt{div}}=2$ and $\epsilon_{\texttt{thr}}=10^{-5}$. 

\begin{table}[h!]
  \begin{center}
    \caption{\gnsindy~performance with different precision value $\epsilon_{\texttt{thr}}$, time division $\texttt{t}_{\texttt{div}}=2$ and fixed sample size $900$ for recovering KdV equation}
    \label{KdV_precision_sensitivity}
    \begin{tabular}{|c|c|c|}
      \toprule 
      \textbf{Precision value} & \textbf{Estimated \pde} \\
      \hline
      \midrule 
      $\epsilon_{\texttt{thr}}=5 \times 10^{-5}$ & $\bu_t  - 5.9994 \bu \bu_x - 0.9995 \bu_{xxx}=0 $ \\
      $\epsilon_{\texttt{thr}}=1\times 10^{-5}$ & $\bu_t  - 5.9955 \bu \bu_x - 0.9975 \bu_{xxx}=0 $ \\
      $\epsilon_{\texttt{thr}}=1 \times 10^{-6}$ & $\bu_t  - 5.9639 \bu \bu_x - 0.9850 \bu_{xxx}=0 $ \\
      $\epsilon_{\texttt{thr}}=1 \times 10^{-7}$&
      $\bu_t  - 5.9860 \bu \bu_x - 0.9959 \bu_{xxx}=0 $ \\
      \bottomrule 
    \end{tabular}
  \end{center}
\end{table}

In a more sophisticated analysis we consider to investigate the impact of  sparse estimator and constraint on the KdV model discovery. To do so, we utilize the knowledge that we already acquired regarding the settings of the \qdeim~algorithm to filter out the informative samples as well as the minimum number of cardinality that is required to recover the model. These settings are $\texttt{t}_{\texttt{div}}=2$, $\epsilon_{\texttt{thr}}=10^{-5}$ and $900$ randomly selected samples out of total $14450$ informative samples. In this regard different combinations of sparse estimator and constraint are considered. The results of these experiments are reported in \Cref{KdV_sensitivity_estimator_constraint} where we see all of the combinations have good performance.

\begin{table}[h!]
  \begin{center}
    \caption{\gnsindy~performance with different sparse estimator and different constraint under fixed \qdeim~setting $\epsilon_{\texttt{thr}}=10^{-5}$, $\texttt{t}_{\texttt{div}}=2$ and fixed sample size $900$ for recovering KdV equation}
    \label{KdV_sensitivity_estimator_constraint}
    \begin{tabular}{|c|c|c|}
      \toprule 
      \textbf{Sparse estimator} & \textbf{Constraint} & \textbf{Estimated \pde}\\
      \hline
      \midrule 
      \stridge & \stridge & $\bu_t  - 5.9955 \bu \bu_x - 0.9975 \bu_{xxx}=0 $\\
      \lasso & \stridge & $\bu_t  - 5.9899 \bu \bu_x - 0.9952 \bu_{xxx}=0 $\\
      \stridge & \ols & $\bu_t  - 5.9579 \bu \bu_x - 0.9827 \bu_{xxx}=0 $
      \\
      \lasso & \ols & $\bu_t  - 5.9799 \bu \bu_x - 0.9933 \bu_{xxx}=0 $\\
      \bottomrule 
    \end{tabular}
  \end{center}
\end{table}


To see how \gnsindy~performs under \dnn~structural variations we consider to do set of experiments. The \gnsindy~hyperparameters are set as mentioned earlier with the only alteration in the number of neurons in each layer. We vary the number of neurons in each layer based on a geometric sequence with a common ratio of $2$ with initial value $8$ neurons. The results of these experiments are reported in \Cref{KdV_NN_sensitivity} where we see \gnsindy~can not recover the KdV when the number of neurons in hidden layers are $8$ and $16$. This reveals the challenge of selecting right \dnn~structure for \pde~discovery.


\begin{table}[h!]
  \begin{center}
    \caption{\gnsindy~performance with different \dnn~structure, fixed precision value $\epsilon_{\texttt{thr}}=10^{-5}$, fixed time division $\texttt{t}_{\texttt{div}}=2$, fixed sample size $900$ for recovering KdV equation, sparse estimator and constraint of type \stridge}
    \label{KdV_NN_sensitivity}
    \begin{tabular}{|c|c|c|}
      \toprule 
      \dnn~structure & \textbf{Estimated \pde} \\
      \hline
      \midrule 
      $[2,\ 8,\ 8,\ 8,\ 8,\ 1]$ & $\bu_t -0.5534 \bu_{xxx} -8.4820 \bu \bu_{x} + 11.1437 \bu^2 \bu_x=0 $
 \\
      $[2,\ 16,\ 16,\ 16,\ 16,\ 1]$ & $\bu_t -1.1290\bu_{xxx}-5.7801 \bu \bu_{x} + 2.3874 \bu \bu_{xxx} -4.7362 \bu^2 \bu_{xxx}=0 $ \\
      $[2,\ 32,\ 32,\ 32,\ 32,\ 1]$ & $\bu_t  - 5.9955 \bu \bu_x - 0.9975 \bu_{xxx}=0 $ \\
      $[2,\ 64,\ 64,\ 64,\ 64,\ 1]$ & $\bu_t  - 5.9812 \bu \bu_x - 0.9923 \bu_{xxx}=0 $ \\
      \bottomrule 
    \end{tabular}
  \end{center}
\end{table}

\begin{figure}[h!]
\centering
\begin{subfigure}{1\textwidth}
    \includegraphics[width=1\textwidth]{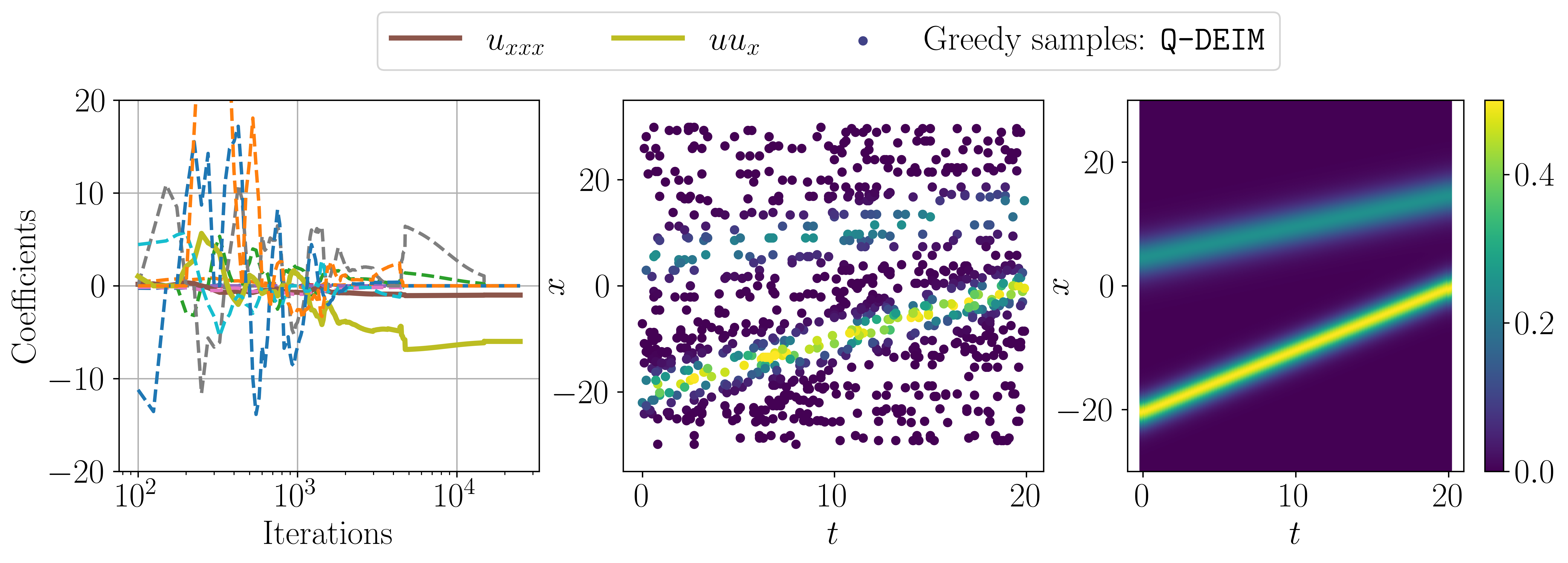}
    \caption{(right) Entire dataset ; (middle) selection of $900$ Greedy samples resulted by \qdeim~algorithm for KdV equation with $\texttt{t}_{\texttt{div}}=2$ and $\epsilon_{\texttt{thr}}=10^{-5}$; (left) estimated coefficients with \gnsindy.}
    \label{fig:KDV_deim_900}
\end{subfigure}
\hfill
\begin{subfigure}{1\textwidth}
    \includegraphics[width=1\textwidth]{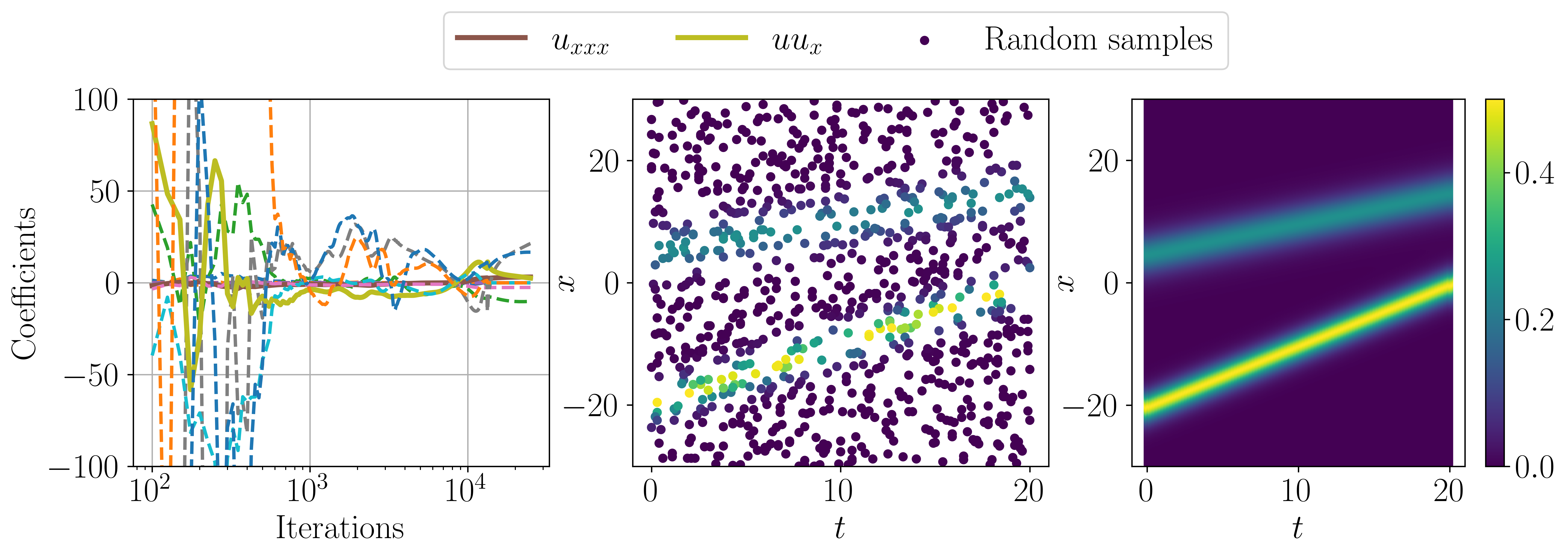}
    \caption{ (right) Entire dataset ; (middle) selection of $900$ random samples for KdV equation; (left) estimated coefficients with \deepymod.
    }
    \label{fig:KDV_random_900}
\end{subfigure}
\caption{ Comparison of the \gnsindy~performance with \deepymod~in KdV equation model discovery}
\end{figure}

Lastly we consider to compare the \gnsindy~performance with \deepymod (\cite{both2021deepmod}). The same \dnn~structure is considered for \deepymod~with its default sparse estimator and constraint of the type \lasso~and \ols~respectively. The results of this simulation is shown in the \Cref{KDV_gnsindy_deepymod} where we clearly see that \gnsindy~outperforms \deepymod~and can recover the \pde~with $900$ greedy samples. It is worth to highlight that \deepymod~selects the samples randomly to train its \dnn~structure. In \Cref{fig:KDV_deim_900} and \Cref{fig:KDV_random_900} the candidate samples that are used in the training loop of each algorithm are shown in the middle graphs, while the left graphs depict the evolution of different coefficients thorough the training loop iterations corresponding to each algorithm. In particular from coefficient evolution of \gnsindy~algorithm we see that the coefficient corresponding to the term $\bu_{xxx}$ has a faster convergence rate respect to the coefficient corresponding to the term $\bu \bu_{x}$.

\begin{table}[h!]
  \begin{center}
    \caption{Comparing \gnsindy~and \deepymod~for KdV \pde~discovery}
    \label{KDV_gnsindy_deepymod}
    \begin{tabular}{|c|c|c|}
      \toprule 
      \textbf{Algorithm} & \textbf{Estimated \pde} \\
      \hline
      \midrule 
     \gnsindy & $\bu_t  - 5.9955 \bu \bu_x - 0.9975 \bu_{xxx}=0 $ \\
      \deepymod \cite{both2021deepmod} & $\bu_t -2.5393\bu_x + 3.1320 \bu_{xxx} + 2.5640\bu \bu_x -10.1184 \bu \bu_{xxx}$ \\ 
      & $+ 21.4131 \bu^2 \bu_x + 16.6569 \bu^2 \bu_{xxx}=0 $
      \\
      \bottomrule 
    \end{tabular}
  \end{center}
\end{table}

\section{Conclusion}\label{sec:conclusion}
A greedy sampling approach has been considered in the framework of sparse identification of non linear dynamical systems (\sindy) which takes advantage of deep neural network (\dnn) and its strength for the \pde~ model discovery. In particular, our proposed methodology is the extension of \deepymod (\cite{both2021deepmod}) with the integration of greedy sampling approach as the data collection and new sparse estimator and constraint functions in the training loops of the \dnn. In this setting, discrete empirical interpolation method (\deim) has been employed to extract the most informative samples of a snapshot matrix associated to a \pde. Due to the usage of greedy samples with the combination of \dnn~and \sindy~algorithm the proposed approach has been named \gnsindy. Our comprehensive study on Burgers' equation, Allen-Cahn equation, and Korteweg-de Vries equation revealed that usage of greedy samples and sequential threshold ridge regression (\stridge) significantly increase the success rate of the model discovery algorithm. In the comparison phase, \gnsindy~outperformed \deepymod~ in all the simulation settings and we could discover Burgers' equation,
Allen-Cahn equation, and Korteweg-de Vries equation with 
with $0.5\%$, $0.1\%$, and $0.874\%$ of the dataset respectively. After conducting a comprehensive analysis with various simulation settings, we have uncovered a fundamental challenge associated with the discovery of \pde~and the selection of appropriate hyperparameters for data-driven approaches. 


\addcontentsline{toc}{section}{References}
\bibliographystyle{IEEEtran}
\bibliography{mybib}

\appendix
\section{Appendix}\label{sec:appendix}

\gnsindy~is implemented as a general framework in \texttt{Python} and is available as a git repository \url{https://gitlab.mpi-magdeburg.mpg.de/forootani/gnsindy} or \url{https://github.com/Ali-Forootani/GN_SINDy}.

\end{document}